\documentclass[12pt,reqno,a4paper]{amsart}

\usepackage{amssymb,latexsym,amsmath,amsthm}
\usepackage{xspace}
\usepackage{rotating}                                                                                                                                                                                                                                                                                                                                                                                                                                                                                                           

\usepackage{epsfig}

\usepackage[french]{babel}

\theoremstyle{plain}
\newtheorem{thm}{Th\'eor\`eme}[section]

\newtheorem{lemma}[thm]{Lemme} 
\newtheorem{corol}[thm]{Corollaire}

\theoremstyle{remark}
\newtheorem{remark}[thm]{Remarque}

\numberwithin{equation}{thm}

\theoremstyle{definition}
\newtheorem{defn}[thm]{D\'efinition}

\usepackage{calc}
\newcounter{hours}
\newcounter{minutes}

\newcommand{\Z}{\mathbb{Z}}

\newcommand{\R}{\mathbb{R}}
\newcommand{\RR}{\mathbb{R}^2}
\newcommand{\difeo}{diff\'eo\-mor\-phis\-me\xspace}
\newcommand{\difeos}{diff\'eo\-mor\-phis\-mes\xspace}

\begin{document}

\title[Localisation des points fixes communs]{Localisation des points fixes communs\\ pour des diff\'eomorphismes commutants du plan}
\author{S. Firmo}

\address{Instituto de Matem\'atica e Estat\'\i stica\\ Universidade Federal Fluminense\\ Br}
\email{firmo@mat.uff.br}

%%%%%%%%%%%%%%%%%%%%%%%%%%%%%%%%%%%%%%%%%%%%%%%%%%%
%%%%%%%%%%%%%%%%%%%%%%%%%%%%%%%%%%%%%%%%%%%%%%%%%%%

%%%%%%%%%%%%%%%%%%%%%%%%%%%%%%%%%%%%%%%%%%%%%%%%%%%
%%%%%%%%%%%%%%%%%%%%%%%%%%%%%%%%%%%%%%%%%%%%%%%%%%%

%\date{Actualis\'e le \today{} \`a \horacerta}
\date{\today{}}



%\thanks{Version pr\'eliminaire}

%%%%%%%%%%%%%%%%%%%%%%%%%%%%%%%%%%%%%%%%%%%%%%%%%%%
%%%%%%%%%%%%%%%%%%%%%%%%%%%%%%%%%%%%%%%%%%%%%%%%%%%

\maketitle

\thispagestyle{empty}

\vskip50pt

%%%%%%%%%%%%%%%%%%%%%%%%%%%%%%%%%%%%%%%%%%%%%%%

%:Resumo
\begin{abstract}
On d\'emontre que si \,$G\subset\text{Diff}^{1}(\RR)$\, est un sous-groupe 
ab\'{e}lien
engendr\'e par une famille quelconque de diff\'eomorphismes de \,$\RR$\, qui sont
\,$C^{1}$-proches de l'identit\'e (pour la \,$C^{1}$-topologie de Whitney) et, 
s'il existe un point \,$p\in\RR$\, dont l'orbite par \,$G$\, est born\'ee, alors les \'{e}l\'{e}ments de 
\,$G$\, ont un point fixe commun dans l'enveloppe convexe de
\,$\overline{\mathcal{O}_{p}(G)}$\,. Ici, \,$\overline{\mathcal{O}_{p}(G)}$\,
est l'adh\'erence topologique de l'orbite de \,$p$\, par \,$G$.
 
\end{abstract}

%%%%%%%%%%%%%%%%%%%%%
\vskip30pt
\section{Introduction}
\vskip10pt
%%%%%%%%%%%%%%%%%%%%%

R\'{e}cemment, Franks-Handel-Parwani  ont d\'emontr\'e dans \cite{fhp01}, que si
\,$G$\,  est un sous-groupe ab\'elien de \,$\text{Diff}^{1}_{+}(\RR)$\,  engendr\'{e} par un nombre fini d'\'el\'ements et, s'il existe un point \,$p\in\RR$\,   dont l'orbite par 
\,$G$\, est born\'ee, alors les \'{e}l\'{e}ments de \,$G$\, poss\`edent un point fixe commun.

En 1989, Bonatti avait d\'emontr\'e l'existence des points fixes communs pour les sous-groupes ab\'{e}liens de \,$\text{Diff}^{1}(S^{2})$\, engendr\'es par des diff\'eomorphismes \,$C^{1}$-proches de l'identit\'e.

En adaptant la strat\'{e}gie de Bonatti dans \cite{bo01}, on donne une version 
du th\'eor\`eme de Franks-Handel-Parwani qui nous permet de localiser le point fixe commun. Cela permet aussi d'avoir une version  valide pour une famille quelconque de \difeos commutants qui sont \,$C^{1}$-proches de l'identit\'e. Pr\'{e}cis\'{e}ment, on d\'emontre le r\'esultat suivant.

%%%%%%%%%%%%%%%%%%%%%%%%%%%%%%%%%
\vskip10pt
\begin{thm}\label{teor:central}
Il existe un voisinage \,$\mathcal{V}$ de l'identit\'e dans 
\,$\mathrm{Diff}^{1}(\RR)$\, muni de la \,$C^{1}$-topologie de Whitney avec la propri\'et\'e suivante.

Si \,$G\subset\mathrm{Diff}^{1}(\RR)$\, est un sous-groupe ab\'elien engendr\'e par une famille quelconque d'\'el\'ements de \,$\mathcal{V}$ et s'il existe un point 
\,$p\in\RR$\, dont l'orbite par \,$G$ est born\'ee, alors les \'el\'ements de 
\,$G$ poss\`edent un point fixe commun dans l'enveloppe convexe de 
\,$\overline{\mathcal{O}_{{p}}(G)}$.

\end{thm}
\vskip5pt
%%%%%%%%%%%%%%%%%%%%%%%%%%%%%%%%%

Dans le cas o\`u le groupe \,$G$ est finiment engendr\'{e}, on peut \^etre un peu plus pr\'ecis par rapport \`a l'ensemble de localisation du point fixe commun 
(voir th\'eor\`eme \ref{teo:cas:fini}).

Le voisinage \,$\mathcal{V}$  cit\'e dans le th\'eor\`eme est un voisinage de type
$$\Big\{f\in\text{Diff}^{1}(\RR) \ ; \ \|f(x)-x\|\,,\|Df(x)-Id\|<\epsilon
\ \text{ sur } \ \RR \Big\}$$
pour \,$\epsilon>0$\, choisi suffisamment petit.

L'id\'{e}e de la preuve est simple. On commence par fixer un ensemble \,$\mathcal{F}$\, de g\'en\'erateurs de \,$G$ form\'e par des 
\difeos \,$C^{1}$-proches de l'identit\'e et on consid\`ere \,$f\in\mathcal{F}$.
Tout d'abord, on garantit l'existence d'un point fixe
\,$q$\, pour \,$f$\, dont l'orbite par les autres g\'en\'erateurs du groupe est contenue dans l'enveloppe convexe de \,$\overline{\mathcal{O}_{p}(G)}$\,. En r\'{e}p\'{e}tant le m\^eme processus avec le point \,$q$\, et un autre \difeo
de \,$\mathcal{F}$, on trouve un point fixe commun pour ces deux \difeos dont l'orbite par les autres g\'en\'erateurs reste dans l'enveloppe convexe de \,$\overline{\mathcal{O}_{p}(G)}$\,. Ainsi, on montre le r\'esultat pour un nombre fini de \difeos de l'ensemble \,$\mathcal{F}$. Le point cl\'e dans ce contexte est d\'emontrer que la \,$G$-orbite du premier point fixe trouv\'e reste contenue dans l'enveloppe convexe de \,$\overline{\mathcal{O}_{p}(G)}$\,. Pour cela, on a besoin de trouver un point fixe pour \,$f$\, dont
l'indice par rapport \`a une famille de courbes  construite \`a  partir de la
\,$f$-orbite du point \,$p$\, soit non nul. 

Ensuite, on passe du cas finiment engendr\'e au cas g\'{e}n\'{e}ral en utilisant la propri\'et\'e d'intersection finie pour une famille de compacts.

Le fait d'\^etre  \,$C^{1}$-proche de l'identit\'e
 est essentielle dans le th\'eor\`eme \ref{teor:central}. On peut construire des \difeos de classe \,$C^{\infty}$ dans le plan, \`a support compact, \,$C^{0}$-proches de l'identit\'e et qui ne satisfont pas le th\'eor\`eme. 

Pour cela, on consid\`ere une application 
\,$\phi_{n}:[0\,,\infty)\rightarrow[0\,,\infty)$\, d\'ecroissante, de classe 
\,$C^{\infty}$ et telle que 
$$\phi_{n}(x)=\begin{cases}
2\pi/n     &  \text{quand} \quad x\in[0\,,1]   \\  
0    &  \text{quand} \quad x\in[2\,,\infty) 
\end{cases}$$
o\`u \,$n\geq2$\, est un entier.
Soit \,$f_{n}$\, le \difeo de \,$\RR$ dont la restriction au cercle
\,$S_{r}:x^{2}+y^{2}=r^{2}$\, est la rotation d'angle \,$\phi_{n}(r)$\, pour 
tout \,$r>0$.
Chaque point dans le disque \,$D_{1}:x^{2}+y^{2}\leq1$\, est un point p\'eriodique de p\'eriode \,$n$\, pour \,$f_{n}$\, sauf l'origine qui est l'unique point fixe sur \,$D_{1}$\,. Pour chaque point \,$p\in D_{1}$\,, l'origine est dans l'enveloppe convexe de l'orbite de \,$p$\, par \,$f_{n}$.

On va modifier \,$f_{n}$\, par conjugaison afin de d\'etruire cette propri\'et\'e. Pour cela, on consid\`{e}re deux points \,$p,q\in\RR$\, 
qui sont \`a une m\^eme distance, positive et plus petite que \,$1$, de l'origine. On suppose aussi que le point \,$q$\, n'appartient pas \`a l'orbite de \,$p$\, 
par \,$f_{n}$. 
Maintenant, on consid\`{e}re un \difeo \,$\psi$\, de \,$\RR$\, \`a support dans un tr\`es petit voisinage \,$\mathcal{U}\subset D_{1}$\, du segment de droite 
 \,$\{tq\in\RR \ ; \ t\in[0\,,1]\}$\, et tel que 
\,$\psi(0)=q$\,. Ainsi, le point \,$q$\, est l'unique point fixe de
 \,$\psi\circ f \circ \psi^{-1}$\,  sur \,$D_{1}$ et, si on a choisi 
\,$\mathcal{U}$\, suffisamment petit, alors
l'orbite de \,$p$\, par 
\,$f$\, et par \,$\psi\circ f \circ \psi^{-1}$\, co\"{\i}ncident. D'autre part, \,$q$\, 
n'est pas dans l'enveloppe convexe de 
l'orbite de \,$p$\, par  \,$\psi\circ f \circ \psi^{-1}$. Bien s\^ur, cette construction peut \^etre r\'ealis\'ee aussi proche qu'on veut de l'identit\'e, dans la \,$C^{0}$-topologie. Il suffit de prendre \,$n$\, grand et \,$p,q$\, proches de l'origine.

%%%%%%%%%%%%%%%%%%%%%
\vskip30pt
\section{Notations}\label{notacoes}
\vskip10pt
%%%%%%%%%%%%%%%%%%%%%

Soient \,$a_{1},\ldots,a_{n}\in\RR$\, des points non n\'{e}cessairement distincts  et tels que \,$a_{i}\neq a_{i+1}$\, pour tout \,$i\in\Z/n\Z$\,.  On note
\,$\Gamma(a_{1},\ldots,a_{n})$\, la courbe ferm\'ee et orient\'ee obtenue en mettant bout \`a bout les segments 
$$[a_{1},a_{2}],[a_{2},a_{3}],\ldots,[a_{n-1},a_{n}],[a_{n},a_{1}]$$
o\`u \,$[a_{i}\,,a_{i+1}]$\, est le segment de droite orient\'e joignant 
\,$a_{i}$\, \`a \,$a_{i+1}$\,. Ces segments seront appel\'es 
\,{\it segments qui composent la courbe}\, \,$\Gamma(a_{1},\ldots,a_{n})$\, 
et chaque \,$a_{i}$\, sera appel\'e  \,{\it sommet}\, de la courbe. On dira aussi que 
\,$[a_{i+1}\,,a_{i+2}]$\, est le \,{\it successeur}\, de \,$[a_{i}\,,a_{i+1}]$\,
dans la courbe, pour tout \,$i\in \Z/n\Z$\,.

On va  supposer aussi que si deux segments orient\'es qui composent la courbe 
\,$\Gamma(a_{1},\ldots,a_{n})$\, s'intersectent, alors ils le font selon un angle inf\'{e}rieur \`a \,$\pi/2$\,. En particulier, s'ils 
s'intersectent selon un segment de droite, alors les orientations co\"{\i}ncident dans l'intersection.
Suivant ces conventions, si 
$$a_{1}=(-1\,,0)=a_{4} \quad ; \quad a_{2}=(1\,,0)=a_{5} \quad ; \quad a_{3}=(0\,,2)=a_{6}$$
donc la courbe \,$\Gamma(a_{1},a_{2},a_{3},a_{4},a_{5},a_{6})$\, est la courbe 
\,$\Gamma(a_{1},a_{2},a_{3})$\, d\'ecrite deux fois dans le sens anti-horaire.
On rappelle que \,$\Gamma(a_{1},a_{2},a_{3},a_{4},a_{5},a_{6})$\, poss\`ede des
auto-intersections qui persistent  par de petites perturbations de la courbe.

\'Etant donn\'e un \difeo \,$f\in \mathrm{Diff}^{1}(\RR)$\, et  un point \,$p\in\RR-\text{Fix}(f)$\,,
on note  \,$\Gamma_{p}^{f}$\, la courbe orient\'ee obtenue en mettant bout \`a bout les segments de la famille  \,$\big\{[f^{n}(p)\,,f^{n+1}(p)]\big\}_{n\in\Z}$\, et
on note  \,$\Gamma_{p,m}^{f}$\, la courbe ferm\'ee, orient\'ee obtenue en mettant bout \`a bout les segments 
$$ [f(p)\,,f^{2}(p)],\ldots,[f^{m-1}(p)\,,f^{m}(p)],[f^{m}(p),f(p)]$$
o\`u \,$m\geq2$\,.

Si \,$\gamma$\, est une courbe param\'etr\'ee, ferm\'ee et continue dans 
\,$\RR$\,  et \,$q\in\RR$\, est un point  hors de la courbe, on notera
\,$\text{Ind}_{q}(\gamma)$\, l'indice de \,$q$\, par rapport \`a \,$\gamma$.
La propri\'et\'e \,$\text{Ind}_{q}(\gamma)\neq0$\, caract\'erise le fait que le point \,$q$\, 
appartient \`a l'enveloppe convexe de l'ensemble des points d\'{e}crits par 
\,$\gamma$.

Nous allons aussi adopter les notations et nomenclatures suivantes:
\begin{enumerate}
\item[--] 
$\mathcal{O}_{p}(f):=\{f^{j}(p)\in\RR \ ; \ j\in\Z\}$;

\noindent
est l'orbite de  \,$p$\, par \,$f$\, ou, simplement, la  \,$f$-orbite de \,$p$:

\item[--]
$\mathcal{O}_{p}(G):=\{g(p)\in\RR  \ ; \ g\in G\}$;

\noindent
est l'orbite de \,$p$\, par \,$G$\, ou, la  \,$G$-orbite de 
\,$p$\, quand \,$G\subset\text{Diff}^{1}(\RR)$\, est un sous-groupe.
 Dans ce cas on utilisera aussi la notation \,$\mathcal{O}_{p}(g_{1},\ldots,g_{m})$\, pour indiquer la  \,$(g_{1},\ldots,g_{m})$-orbite de 
\,$p$\, quand \,$\{g_{1},\ldots,g_{m}\}$\, est un ensemble de g\'en\'erateurs 
de \,$G$;

\item[--]
si \,$A\subset\RR$, on note \,$\text{Conv}(A)$\, (resp.
\,$\overline{A}$\,)\, l'enveloppe convexe (resp. l'adh\'erence) de 
\,$A$;

\item[--]
on note \,$\text{Fix}(f)$\,, \,$\text{Fix}(f_{1},\ldots,f_{k})$\, et 
\,$\text{Fix}(\mathcal{F},\mathcal{G})$\, les ensembles des points fixes communs de
\,$f$\,, de 
\,$f_{1},\ldots,f_{k}$\, et de \,$\mathcal{F} \cup \mathcal{G}$, respectivement,
o\`u \,$\mathcal{F},\mathcal{G}\subset \text{Diff}^{1}(\RR)$.

\end{enumerate}

Dans \,$\RR$ on utilise la norme usuelle qui sera not\'ee 
\,$\|{\cdot}\|$\,. L'ensemble
\,$\text{Diff}^{1}(\RR)$\, est l'espace des \difeos de 
\,$\RR$\, de classe \,$C^{1}$ munie de la \,$C^{1}$-topologie de Whitney.
On notera \,$B(q\,;\delta)$\, la boule ferm\'ee de centre  \,$q\in\RR$\, et de rayon \,$\delta>0$\, pour la norme \,$\| {\cdot} \|$\,.

%%%%%%%%%%%%%%%%%%%%%
\vskip30pt
\section{Un lemme topologique}\label{lema:top}
\vskip10pt
%%%%%%%%%%%%%%%%%%%%%

Dans cette section, nous  \'{e}non\c cons un lemme topologique
qui aura un r\^{o}le important pour la preuve du th\'eor\`eme  \ref{teor:central}. 
La preuve du lemme sera pr\'esent\'ee dans la derni\`ere section de cet article.

Dans la suite, on suppose que \,$a_{1},\ldots,a_{n}\in\RR$\, sont des points non n\'{e}cessairement distincts mais tels que \,$a_{i}\neq a_{i+1}$\, quand
\,$i\in \Z/n\Z$\, o\`u \,$n\geq2$\,. On admet aussi que l'angle entre deux 
segments \,$[a_{i},a_{i+1}]$\, et \,$[a_{j},a_{j+1}]$\, qui s'intersectent est 
inf\'erieur \`a \,$\pi/2$\, pour tout \,$i\,,j\in \Z/n\Z$\,.

%%%%%%%%%%%%%%%%%%%%%%%%%%%%%%%%%%
%:Lema Topol—gico {deco:curv:circ}
\vskip10pt
\begin{lemma}\label{deco:curv:circ}
La courbe \,$\Gamma=\Gamma(a_{1},\ldots,a_{n})$\, se d\'ecompose en un nombre fini de courbes simples ferm\'ees
\,$\gamma_{1},\ldots,\gamma_{m}$\, ayant les propri\'et\'es$\,:$
\begin{enumerate}
\item[$(i)$]
le nombre de fois que chaque point du plan est recouvert par \,$\Gamma$\,
est le m\^eme que par \,$\gamma_{1}\cup\cdots\cup\gamma_{m}\,;$
\item[$(ii)$]
$\gamma_{j}=\Gamma(b_{j,1},\ldots,b_{j,n_{j}})$\, o\`u chaque segment  qui compose  \,$\gamma_{j}$\, est contenu dans un des segments qui compose \,$\Gamma$,
ayant la m\^eme orientation que celui-ci, pour chaque \,$j\in\{1,\ldots,m\}$
$($dans ce sens on dira que \,$\gamma_{j}\subset \Gamma)\,;$

\item[$(iii)$]
si les int\'{e}rieurs des disques bord\'es par \,$\gamma_{i}$\, et \,$\gamma_{j}$\,
respectivement, s'intersectent, alors l'un d'entre eux est contenu 
dans l'autre$\,;$

\item[$(iv)$]
$\displaystyle\mathrm{Ind}_{q}(\Gamma)
=\sum_{1\leq j\leq m}
\mathrm{Ind}_{q}(\gamma_{j})$\, pour tout
\,$q\in\RR-\Gamma\,;$

\item[$(v)$]
il existe \,$\kappa\in\{1,\ldots,m\}$\, tel que \,$\mathrm{Ind}_{q}(\Gamma)$\, 
est bien d\'efini et  non nul pour tout point \,$q$\, dans l'int\'erieur du disque bord\'e par \,$\gamma_{\kappa}$\,.

\end{enumerate}

\end{lemma}
\vskip5pt
%%%%%%%%%%%%%%%%%%%%%%%%%%%%%%%%%%

Comme corollaire imm\'ediat nous avons le r\'{e}sultat suivant.

%%%%%%%%%%%%%%%%%%%%%%%%%%%%%%%%%%
\vskip10pt
\begin{corol}\label{corol:deco:curv:circ}
\'Etant donn\'e \,$\Gamma=\Gamma(a_{1},\ldots,a_{n})$\,, il existe une courbe simple ferm\'ee \,$\gamma\subset\Gamma$\, telle que 
\,$\mathrm{Ind}_{q}(\Gamma)$\, est bien d\'efini et non nul pour tout point 
\,$q$\, dans l'int\'erieur du disque bord\'e par \,$\gamma$.

\end{corol}
\vskip5pt
%%%%%%%%%%%%%%%%%%%%%%%%%%%%%%%%%%

Afin de pr\'eparer la preuve du lemme topologique nous pr\'esentons ici un m\'ethode de modification d'une courbe \,$\Gamma(a_{1},\ldots,a_{n})$.

Nous allons d\'ecomposer la courbe \,$\Gamma=\Gamma(a_{1},\ldots,a_{n})$\,
en d\'{e}coupant de fa\c con convenable les segments qui composent la courbe, en des segments plus petits. Nous autorisons aussi \`a modifier l'ordre selon laquelle ces segments sont d\'ecrits. Par contre, au cours des modifications, le support de la courbe, la multiplicit\'e de ses points et l'orientation des segments qui composent la courbe seront  pr\'eserv\'es.
 Les modifications seront faites pour d\'etruire les auto-intersections persistantes par de petites perturbations de la courbe  \,$\Gamma$. En plus, dans chaque pas du processus de modification de la courbe, nous serons toujours \`a produire des courbes du type
\,$\Gamma(b_{1},\ldots,b_{s})$\, contenues dans \,$\Gamma$\,.

Regardons le cas particulier o\`u \,$\Gamma$\, poss\`ede seulement un nombre fini de points multiples. Dans ce cas, en augmentant de fa\c con convenable le nombre de sommets de la courbe,  on peut admettre sans perte de g\'{e}n\'{e}ralit\'{e} que si deux segments  qui composent la courbe s'intersectent, alors ils le font selon une extr\'{e}mit\'{e} commune. Dans ce cas,
chaque point de la courbe \,$\Gamma$\, est recouvert une seule fois, sauf les sommets qui peuvent l'\^etre  plusieurs fois.

Prenons alors un point \,$b$\, de la courbe \,$\Gamma$\, o\`u on a des intersections. Ces intersections peuvent \^etre persistantes ou pas par de petites perturbations. Rappelons que par hypoth\`ese, si deux segments qui composent la courbe s'intersectent, alors l'angle entre eux est plus petit que \,$\pi/2$\,.
Ainsi, on peut fixer une ligne \,$L$\, passant par \,$b$\, convenablement orient\'ee, de telle fa\c con que tous les segments qui composent la courbe et qui pointent 
vers \,$b$\, 
$$[\xi_{1},b],\ldots,[\xi_{k},b]$$
sont du c\^ot\'e  droit de \,$L$\, et les segments qui sortent de \,$b$\,
$$[b\,,\nu_{1}],\ldots,[b\,,\nu_{k}]$$ 
sont \`a gauche de \,$L$ comme  montr\'e dans la figure \`a suivre.
\'{E}videmment, ces deux classes de segments ont le m\^eme nombre d'\'el\'ements et
chaque segment  \,$[\xi_{i},b]$\, a comme successeur, dans la courbe \,$\Gamma$,
un segment \,$[b,\nu_{j}]$\,. Maintenant, il est facile de d\'{e}truire, simultan\'{e}ment, tous les intersections persistantes, qui
arrivent au point \,$b$: il suffit de modifier la courbe en mettant \,$[b\,,\nu_{i}]$\,
comme successeur de \,$[\xi_{i}\,,b]$\,. Bien s\^ur, on suppose que \,$\nu_{j}$\, 
est \`a droite de  \,$\nu_{j+1}$\, et que 
\,$\xi_{i}$\, est \`a droite de \,$\xi_{i+1}$\,. 
Apr\`es cette modification, tous les intersections au point \,$b$\, peuvent 
\^etre d\'etruites  par de petites perturbations de la courbe.

\vskip10pt

\setlength{\unitlength}{1mm}%
\hfil
%\fbox{%
\begin{picture}(32,38)(-1,-4)%

\put(0,0){\includegraphics[width=3cm]{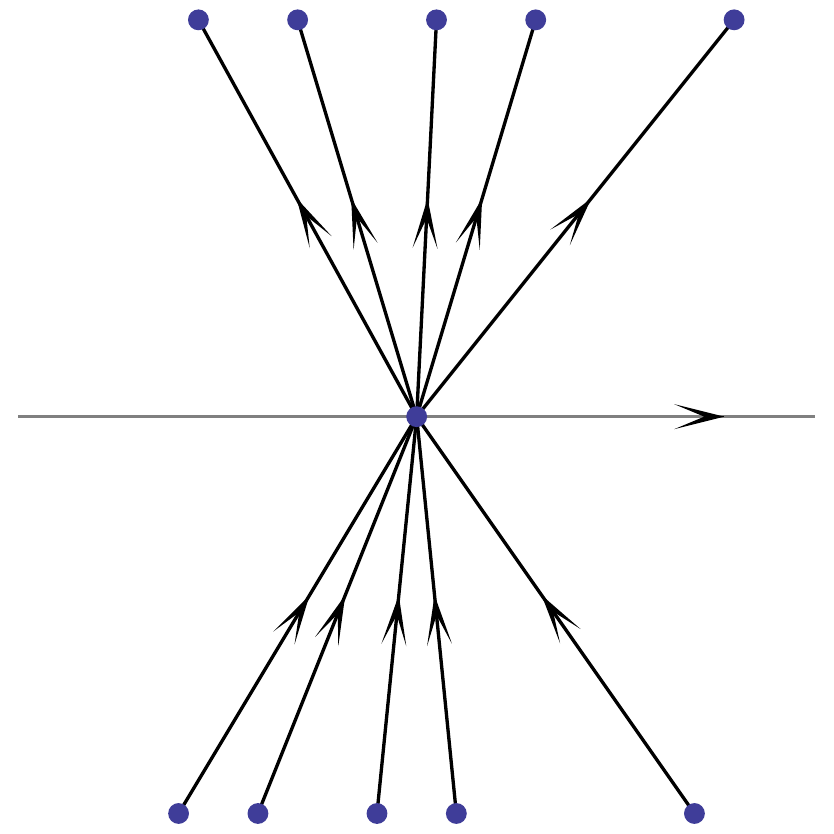}}

\put(12,14){\scriptsize$b$}
\put(28,16){\scriptsize$L$}
\put(24,-2.5){\tiny$\xi_{1}$}
\put(16,-2.5){\tiny$\xi_{2}$}
\put(12,-2.5){\tiny$\xi_{3}$}
\put(4,-2.5){\tiny$\xi_{k}$}

\put(8,-2){\tiny$\ldots$}

\put(25,30.5){\tiny$\nu_{1}$}
\put(19,30.5){\tiny$\nu_{2}$}
\put(14,30.5){\tiny$\nu_{3}$}
\put(4,30.5){\tiny$\nu_{k}$}

\put(9,31){\tiny$\ldots$}

\end{picture}%
%}

\vskip10pt

La condition $(iii)$ sera essentielle pour d\'emontrer $(v)$ comme 
les figures \`a suivre le montrent. Dans la figure \`a gauche, on a la courbe 
\,$\Gamma(a_{1},\ldots,a_{10})$. Apr\`es faire une premi\`ere modification, on obtient les courbes \,$\Gamma(a_{6},b_{1},a_{10},a_{1},\ldots,a_{5})$\, et 
\,$\Gamma(a_{9},b_{1},a_{7},a_{8})$\, qui sont montr\'ees dans la figure au milieu.
Les deux courbes obtenues sont des courbes simples  ferm\'ees, mais elles ne satisfont pas la condition $(iii)$. En fait, dans l'intersection des int\'{e}rieurs des disques bord\'es par ces courbes, chaque point a un indice nul par rapport \`a la courbe 
\,$\Gamma(a_{1},\ldots,a_{10})$\,.

Mais, il reste encore deux intersections persistantes par de petites perturbations:
$$[a_{3}\,,a_{4}]\cap[a_{8}\,,a_{9}] \quad \text{et} \quad  
[a_{4}\,,a_{5}]\cap[a_{9}\,,b_{1}] \,.$$
En modifiant les courbes dans ces intersections pour \'{e}liminer les intersections persistantes, on obtient les courbes suivantes:
$$\Gamma(a_{6},b_{1},a_{10},a_{1},a_{2},a_{3},b_{2},a_{9},b_{3},a_{5}) 
\quad \text{et} \quad \Gamma(b_{1},a_{7},a_{8},b_{2},a_{4},b_{3})$$
qui sont montr\'ees dans la figure \`a droite. Ces courbes satisfont les conditions du lemme. Les points dans les int\'{e}rieurs
des disques bord\'es par ces deux courbes  ont des indices non nuls par rapport 
\`a la courbe initiale \,$\Gamma(a_{1},\ldots,a_{10})$\,.

%\vskip10pt

\setlength{\unitlength}{1mm}%
\noindent
\hfil
%\fbox{%
\begin{picture}(45,34)(-1,-1)%

\put(0,0){\includegraphics[width=4cm]{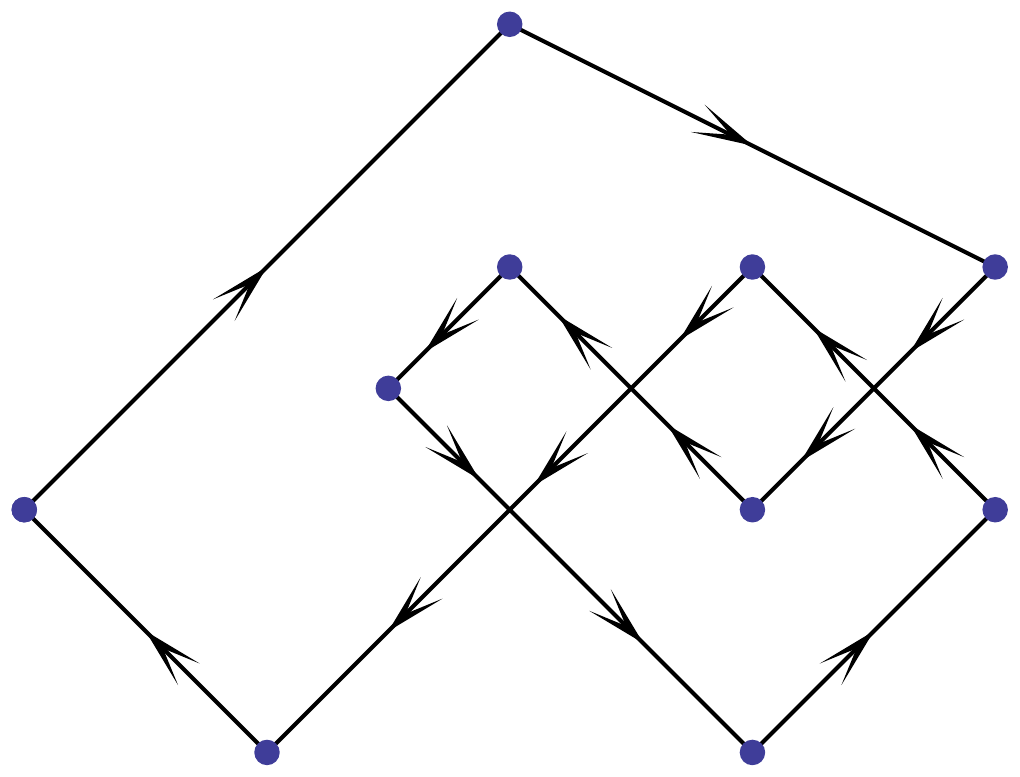}}

\put(2.5,9.7){\scriptsize$a_{1}$}
\put(19,30.5){\scriptsize$a_{2}$}
\put(40,19){\scriptsize$a_{3}$}
\put(30,9){\scriptsize$a_{4}$}
\put(19,21){\scriptsize$a_{5}$}
\put(12,13){\scriptsize$a_{6}$}
\put(30,-1){\scriptsize$a_{7}$}
\put(40,9){\scriptsize$a_{8}$}
\put(29,21){\scriptsize$a_{9}$}
\put(11,-1){\scriptsize$a_{10}$}

\end{picture}%
%}
\setlength{\unitlength}{1mm}%
\hfil
%\fbox{%
\begin{picture}(45,34)(-1,-1)%

\put(0,0){\includegraphics[width=4cm]{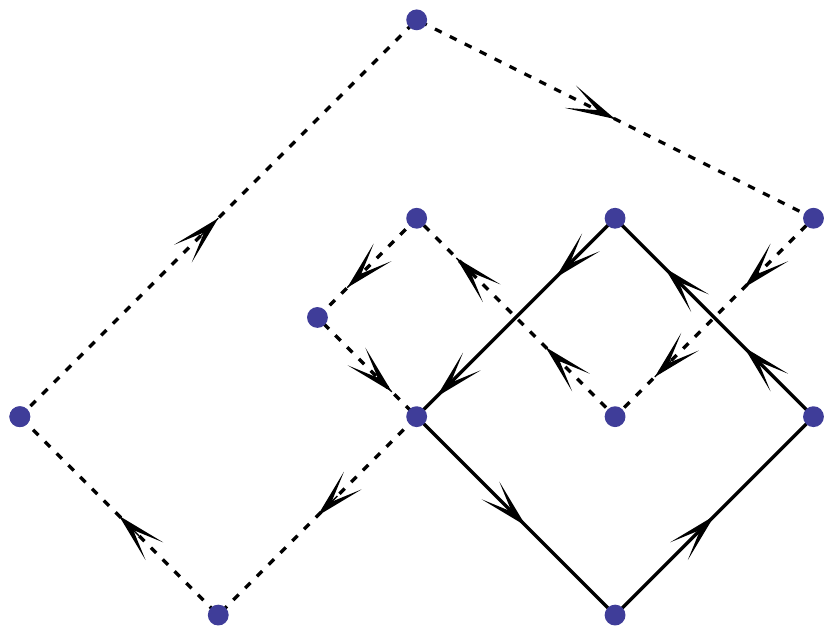}}

\put(2.5,9.7){\scriptsize$a_{1}$}
\put(19,30.5){\scriptsize$a_{2}$}
\put(40,19){\scriptsize$a_{3}$}
\put(30,9){\scriptsize$a_{4}$}
\put(19,21){\scriptsize$a_{5}$}
\put(12,13){\scriptsize$a_{6}$}
\put(30,-1){\scriptsize$a_{7}$}
\put(40,9){\scriptsize$a_{8}$}
\put(29,21){\scriptsize$a_{9}$}
\put(11,-1){\scriptsize$a_{10}$}
\put(19,7){\scriptsize$b_{1}$}

\end{picture}%
%}
\setlength{\unitlength}{1mm}%
\hfil
%\fbox{%
\begin{picture}(45,34)(-1,-1)%

\put(0,0){\includegraphics[width=4cm]{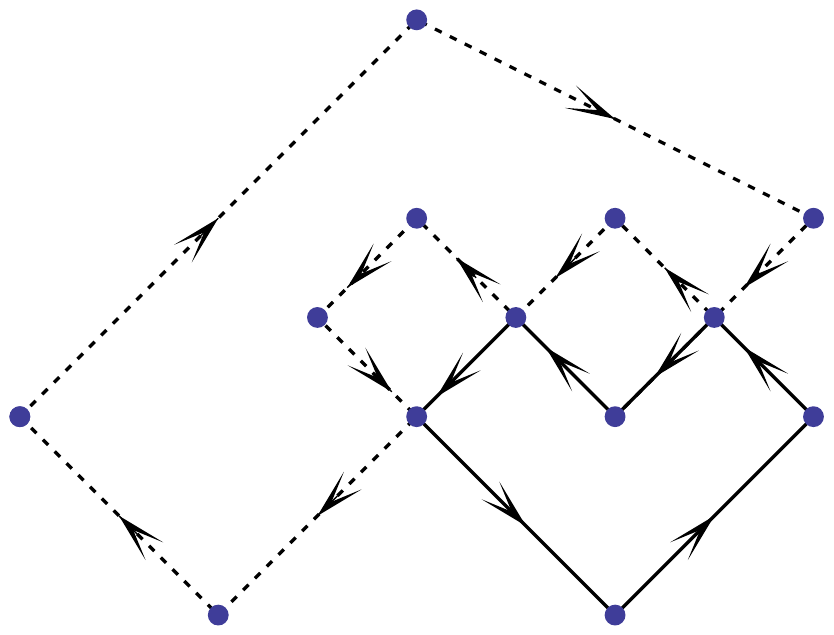}}

\put(2.5,9.7){\scriptsize$a_{1}$}
\put(19,30.5){\scriptsize$a_{2}$}
\put(40,19){\scriptsize$a_{3}$}
\put(30,9){\scriptsize$a_{4}$}
\put(19,21){\scriptsize$a_{5}$}
\put(12,13){\scriptsize$a_{6}$}
\put(30,-1){\scriptsize$a_{7}$}
\put(40,9){\scriptsize$a_{8}$}
\put(29,21){\scriptsize$a_{9}$}
\put(11,-1){\scriptsize$a_{10}$}
\put(19,7){\scriptsize$b_{1}$}
\put(36,14.3){\scriptsize$b_{2}$}
\put(26.5,14.3){\scriptsize$b_{3}$}

\end{picture}%
%}
\hfil\null

\vskip10pt

La d\'{e}composition, dans le cas de la courbe 
\,$\Gamma(a_{1},a_{2},a_{3},a_{4},a_{5},a_{6})$\,
d\'{e}finie dans la section \ref{notacoes} et o\`u 
$$a_{1}=(-1\,,0)=a_{4} \quad ; \quad a_{2}=(1\,,0)=a_{5}\quad \text{e} \quad a_{3}=(0\,,2)=a_{6}\,,$$ 
est donn\'ee par les courbes 
$$\beta_{1}=\Gamma(a_{1},a_{2},a_{3}) \quad \text{et} \quad 
\beta_{2}=\Gamma(a_{4},a_{5},a_{6}) \,.$$

%%%%%%%%%%%%%%%%%%%%%
\vskip30pt
\section{Propri\'et\'es des diff\'eomorphismes \,$C^{1}$-proches de l'identit\'e
}\label{sec:prop:difeo}
\vskip10pt
%%%%%%%%%%%%%%%%%%%%%

Dans cette section, nous \'enon\c cons quelques-uns des r\'esultats de Bonatti (voir \cite{bo01}), ici \'enonc\'es  pour les \difeos de \,$\RR$\, qui sont 
\,$C^{1}$-proches de l'identit\'e dans la \,$C^{1}$-topologie de Whitney.

Dans la suite, \,$\mathfrak{C}^{1}_{s}(\R^{2})$\, est l'espace des applications
\,$f:\R^{2}\rightarrow\R^{2}$\, de classe \,$C^{1}$ muni de la 
\,$C^{1}$-topologie de Whitney. Pour chaque \,$\epsilon>0$\,, on notera 
\,$\mathcal{U}_{\epsilon}\subset \mathfrak{C}^{1}_{s}(\R^{2})$\, le 
voisinage de l'identit\'e donn\'e par 
$$\mathcal{U}_{\epsilon}:=\Big\{f\in \mathfrak{C}^{1}_{s}(\R^{2}) \ ; \ 
\max\big\{\|f(x)-x\|\,,\|Df(x)-Id\|\big\}<\epsilon \ , \ \forall x\in\R^{2} \Big\}$$
o\`u \,$\|Df(x)-Id\|=\sup\big\{\|Df(x)\cdot v-v\|  \ ; \ \|v\|=1
\big\}$\,.

\vskip5pt

Si \,$\epsilon>0$\, est choisi suffisamment petit, \,$\mathcal{U}_{\epsilon}$\,
s'inclut dans 
\,$\text{Diff}^{1}(\R^{2})$\,. Dans la suite, on supposera toujours  cette condition 
v\'erifi\'ee.

\vskip10pt

{\it Dans l'article cit\'e ci-dessus, Bonatti montre qu'il existe \,$\epsilon_{1}>0$\,
avec les propri\'et\'es suivantes}:
\begin{enumerate}
\vskip5pt

\item[$(B.1)$]
{\it Si \,$f\in \mathcal{U}_{\epsilon_{1}}$\, et \,$p\in \RR-\text{\rm Fix}(f)$\,, alors
\,$f$\, est sans point fixe sur la boule ferm\'ee de centre 
\,$p$ et de rayon \,$4\,\|f(p)-p\|$};

\vskip5pt

\item[$(B.2)$]
{\it Si \,$f\in \mathcal{U}_{\epsilon_{1}}$\, et \,$p\in \RR-\text{\rm Fix}(f)$\, alors l'angle entre deux segments orient\'es
$$[f^{i}(p)\,,f^{i+1}(p)] \quad \text{\rm et} \quad  
[f^{j}(p)\,,f^{j+1}(p)]$$ 
qui s'intersectent est plus petit que 
\,$\pi/4$\, pour tout \,$i\,,j\in\Z$};

\vskip5pt

\item[$(B.3)$]
{\it Soient \,$g_{1},\ldots,g_{m}\,,f\in\mathcal{U}_{\epsilon_{1}}$ des \difeos
commutants et soit \,$p\in\RR$\, un point fixe commun pour \,$g_{1},\ldots,g_{m}$\, tel que 
\,$f(p)\neq p$. Supposons qu'il existe une courbe simple ferm\'ee 
 \,$\gamma$\, contenue dans \,$\Gamma^{f}_{p}$\,. Alors, il existe un point fixe commun \`a \,$g_{1},\ldots,g_{m}\,,f$\, dans l'int\'{e}rieur du disque bord\'e 
 par \,$\gamma$.}

\end{enumerate}

\vskip10pt

Les r\'esultats \'enonc\'es dans $(B.1)\,,(B.2)$ et $(B.3)$ correspondent, 
respectivement, aux corollaire 2.2, lemme 2.3 et lemme 5.1, d\'emontr\'es
dans  \cite{bo01}.

On rappelle que la courbe \,$\gamma$\, (dans l'assertion B.3) existe toujours si \,$p$\, est un point 
\,$\omega$-r\'ecurrent pour \,$f$, c'est-\`a-dire, dans ce cas, la courbe 
\,$\Gamma^{f}_{p}$\, poss\`ede au moins un point multiple. Cela a d\'ej\`a  \'et\'e remarqu\'e par Bonatti.

Le prochain lemme n'est pas explicitement d\'emontr\'e dans \cite{bo01}, mais
sa preuve est une r\'{e}p\'{e}tition exacte des raisonnements qui am\`enent \`a la preuve
du lemme 5.1 dans  \cite{bo01}.

%%%%%%%%%%%%%%%%%%%%%%%%%%%%%%%%%%
%:Lema de Bonatti mais geral {bona:mais:geral:01}
\vskip10pt
\begin{lemma}\label{bona:mais:geral:01}
Soient  \,$g_{1},\ldots,g_{m}\,,f\in\mathcal{U}_{\epsilon_{1}}$\, des \difeos
commutants et soit \,$p\in\RR$\, un point fixe commun \`a 
\,$g_{1},\ldots,g_{m}$\, tel que \,$f(p)\neq p$.
Supposons qu'il existe une suite croissante \,$(n_{k})_{k\geq1}$\, dans 
\,$\Z^{+}$\, telle que \,$f^{n_{k}}(p)\rightarrow p$\, et soit \,$\gamma$\, une courbe simple ferm\'ee contenue dans \,$\Gamma^{f}_{p,n_{k}}$\,.

Alors, si  \,$k$\, 
est suffisamment grand, il existe un point fixe commun \`a \,$g_{1},\ldots,g_{m}\,,f$\,  dans l'int\'{e}rieur du disque bord\'e par \,$\gamma$.
\end{lemma}
\vskip5pt
%%%%%%%%%%%%%%%%%%%%%%%%%%%%%%%%%%

\begin{proof}[Commentaire sur la preuve]
La preuve du lemme \ref{bona:mais:geral:01} est une adaptation simple de la preuve du lemme 5.1 de
 \cite{bo01}. La diff\'{e}rence est que dans notre lemme, la courbe simple ferm\'ee \,$\gamma$\, peut avoir des segments de droite contenus dans le segment  \,$[f^{n_{k}}(p),f(p)]$\,. Mais \c ca n'impose pas de 
restrictions car, comme l'a d\'emontr\'e Bonatti dans \cite{bo01}, les propri\'et\'es de \,$f$\, sur les segments \,$[f^{n_{k}}(p),f(p)]$\,, utilis\'{e}es dans la preuve, sont similaires aux propri\'et\'es de \,$f$\, sur le segment \,$[p,f(p)]$\,
quand \,$k$\, est suffisamment grand. En fait, dans \cite{bo02}, Bonatti utilise 
d\'ej\`a les courbes \,$\Gamma^{f}_{p,n_{k}}$\, pour d\'emontrer l'existence des points fixes communs pour des \difeos commutants des surfaces compactes \`a caract\'{e}ristique d'Euler non nulle.
\end{proof}

En combinant les lemmes \ref{deco:curv:circ} et \ref{bona:mais:geral:01} on obtient le r\'{e}sultat suivant.

%%%%%%%%%%%%%%%%%%%%%%%%%%%%%%%%%%
%:Lema de Bonatti com indice n‹o nulo {lema:fundam:ind}
\vskip10pt
\begin{lemma}\label{lema:fundam:ind}
Soient  \,$g_{1},\ldots,g_{m}\,,f\in\mathcal{U}_{\epsilon_{1}}$\, des \difeos
commutants et soit \,$p\in\RR$\, un point fixe commun \`a
\,$g_{1},\ldots,g_{m}$\, tel que \,$f(p)\neq p$. 
Supposons qu'il existe une suite croissante \,$(n_{k})_{k\geq1}$\, dans 
\,$\Z^{+}$ tel que \,$f^{n_{k}}(p)\rightarrow p$.
Alors, si  \,$k$\, 
est suffisamment grand, il existe un point \,$q_{k}\in \mathrm{Fix}(g_{1},\ldots,g_{m}\,;f)$\, tel que 
\,$\mathrm{Ind}_{q_{k}}(\Gamma^{f}_{p,n_{k}})\neq 0$\,.

\end{lemma}
\vskip5pt
%%%%%%%%%%%%%%%%%%%%%%%%%%%%%%%%%%

\begin{proof}[Preuve]
Soient \,$p\in\text{Fix}(g_{1},\ldots,g_{m})-\text{Fix}(f)$\, et 
\,$(n_{k})_{k\geq1}$\,  v\'erifiant les conditions du lemme.
Pour chaque 
\,$k\geq1$\,  consid\'{e}rons la courbe  \,$\Gamma^{f}_{p,n_{k}}$\,.

Rappelons que le segment \,$[f^{n_{k}}(p),f(p)]$\, s'approche de \,$[p,f(p)]$\,
si on prend \,$k$\,  grand. Ainsi d'apr\`es $(B.2)$, l'angle entre deux segments qui composent la courbe \,$\Gamma^{f}_{p,n_{k}}$\,
et qui s'intersectent est inf\'erieur \`a \,$\pi/3$\, si \,$k$\, est suffisamment grand. Alors, le corollaire \ref{corol:deco:curv:circ} nous garantit qu'il existe une courbe simple ferm\'ee \,$\gamma$\, 
contenue dans \,$\Gamma^{f}_{p,n_{k}}$\,. En plus, l'indice 
\,$\text{Ind}_{q}(\Gamma^{f}_{p,n_{k}})$\,
est bien d\'{e}fini   et  non nul pour tous les points dans l'int\'{e}rieur du disque \,$D$\, bord\'e par \,$\gamma$. D'autre part, le lemme \ref{bona:mais:geral:01}
nous garantit qu'il existe un point fixe commun 
\,$q_{k}\in\text{Fix}(g_{1},\ldots,g_{m}\,;f)$\, dans l'int\'{e}rieur de \,$D$\, et
la preuve est termin\'ee.
\end{proof}

Le voisinage \,$\mathcal{V}$\, du th\'eor\`eme \ref{teor:central} ne sera pas le voisinage \,$\mathcal{U}_{\epsilon_{1}}$\, utilis\'{e} dans les r\'esultats de 
\cite{bo01}. Il sera modifi\'e par le voisinage d\'{e}crit dans le prochain lemme.

%%%%%%%%%%%%%%%%%%%%%%%%%%%%%%%%%%
%:Lema tecnico sobre C1-difeomorfismo: prop:loc:dif
\vskip10pt
\begin{lemma}\label{prop:loc:dif}
\'Etant donn\'e  \,$\delta>0$\, et \,$K>1$\,, il existe
\,$\epsilon_{2}=\epsilon_{2}(\delta,K)>0$\, avec la propri\'et\'e suivante$\,:$
$$\|f(\lambda)-f(p)\|\leq K \|f(q)-f(p)\| $$
pour tout \,$f\in\mathcal{U}_{\epsilon_{2}}$\,, 
\,$\|p-q\|\leq \delta$\, et \,$\lambda \in [p\,,q]$\,.
\end{lemma}
\vskip5pt
%%%%%%%%%%%%%%%%%%%%%%%%%%%%%%%%%%

Le lemme \ref{prop:loc:dif} est une cons\'{e}quence imm\'{e}diate du 
lemme~\ref {loc:impt:01}. Dans la suite, pour 
\,$f:{B(0\,,r)}\subset\R^{2}\rightarrow \R^{2}$\,  de classe \,$C^{1}$, la propri\'et\'e  \,$\|f-Id\|_{1}<\epsilon_{2}$\, sur 
\,${B(0\,,r)}\subset\R^{2}$\, signifie que 
$$\sup\Big\{ \|f(x)-x\| \,, \|Df(x)-Id\| \ ; \ x\in {B(0\,,r)}
\Big\}<\epsilon_{2} \,.$$

%%%%%%%%%%%%%%%%%%%%%%%%%%%%%%%%%%
\vskip10pt
\begin{lemma}\label{loc:impt:01}
\'Etant donn\'e \,$r>0$\, et  \,$K>1$\,, il existe 
\,$\epsilon_{2}=\epsilon_{2}(r,K)>0$\, avec la propri\'et\'e suivante$\,:$
si \,$f:{B(0\,,r)}\subset\R^{2}\rightarrow \R^{2}$\, est 
 de classe \,$C^1$  et 
\,$\|\,f-Id\,\|_{1} < \epsilon_{2}$\, sur \,${B(0\,,r)}$\,, alors
pour tous \,$p,q \in {B(0\,,r)}$\, et pour tout 
\,$\lambda\in[\,p\,,q\,]$\, on a
$$\|f(\lambda)-f(p)\|
 \leq K \, \|f(q)-f(p)\| \,.$$
\end{lemma}
\vskip5pt
%%%%%%%%%%%%%%%%%%%%%%%%%%%%%%%%%%

\begin{proof}[Preuve]
On fera une preuve par contradiction. Supposons que le lemme n'est pas vrai. Dans ce cas, pour chaque \,$n\in\mathbb{Z}^{+}$, il existe: 
\begin{enumerate}
\item[--] 
$f_n:{B(0\,,r)}\rightarrow\RR$\, de classe \,$C^{1}$ 
 tel que \,$\|\,f_n-Id\,\|_1 < 1/n$\, 
sur \,${B(0\,,r)} \,;$

\item[--]
$p_n\,,q_n \in {B(0\,,r)}$\, distincts et 
\,$\lambda_n \in [p_n \,, q_n]$\, avec \,$\lambda_{n}\neq p_{n},q_{n} \,;$
\end{enumerate}
tels que
$$ \|f_n(\lambda_n)-f_n(p_n)\|
 > K \, \|f_n(q_n)-f_n(p_n)\|$$
et par cons\'equent,
\begin{equation}\label{desig:lem:loc}
\begin{split} & \frac{1}{\|\lambda_n- p_n\|} \,\big\|f_n (\lambda_n)
   -f_n(p_n) \big\|  >
 \frac{K}{\|q_n-p_n\|} \,
\big\|f_n(q_{n})-f_n(p_n)\big\|
%\notag 
\end{split}
\end{equation}

\vskip5pt

\noindent 
pour tout \,$n\in\mathbb{Z}^+$.
De plus, quitte \`a consid\'erer une sous-suite,  on peut supposer que
\,$(p_n)_{n\geq1}\,,(q_n)_{n\geq1}$\, et \,$(\lambda_{n})_{n\geq1}$\, convergent respectivement vers 
\,$p,q$\, et \,$\lambda$\, 
et que 
$$\bigg(\frac{\lambda_n -p_n}{ \|\lambda_n - p_n\|}\bigg)_{n\geq1} \quad \text{et} \quad 
\bigg(\frac{q_n -p_n}{ \|q_n - p_n\|}\bigg)_{n\geq1}$$ 
convergent, simultan\'{e}ment, vers un m\^eme vecteur unitaire \,$v\in\mathbb{R}^2$.

Si \,$\|q-p\|$\, et \,$\|\lambda-p\|$\, sont strictement positifs alors, en revenant \`a
l'in\'{e}galit\'{e} \eqref{desig:lem:loc} et en passant \`a la limite pour 
\,$n\rightarrow\infty$ on  d\'eduit que l'on a
\begin{align}\label{contra:01}
 1=\| v\|\geq K \, \| v\| = K \,,
 \end{align} 
ce qui nous donne la contradiction cherch\'ee pour ce cas.

Supposons  que \,$\|\lambda_{n}-p_{n}\|\rightarrow 0$\, et \,$\|q-p\|>0$. Par le 
\,{\it th\'eor\`eme des accroissements finis}\, appliqu\'e aux coordonn\'{e}es
\,$f_{n,1}\,,f_{n,2}$\, de \,$f_{n}$\,, on conclut qu'il existe des points 
\,$\xi_{n,1},\xi_{n,2}$\, dans le segment \,$[p_{n}\,,\lambda_{n}]$\,
tels que
\begin{equation}\label{desig:lem:loc:02}
\begin{split} & \frac{1}{\|\lambda_n- p_n\|} \,\Big\{f_n(\lambda_{n})
   -f_n(p_n)\Big\} =  \\  & \qquad =
\bigg(Df_{n,1}(\xi_{n,1})\cdot \frac{\lambda_{n}-p_{n}}{\|\lambda_n- p_n\|}
\, , \, Df_{n,2}(\xi_{n,2})\cdot \frac{\lambda_{n}-p_{n}}{\|\lambda_n- p_n\|}
\bigg)\,.
%\notag 
\end{split}
\end{equation}
Comme  
\,$\|\lambda_{n}-p_{n}\|\rightarrow 0$\, et \,$p_{n}\rightarrow p$\,,
il s'en suit que \,$\xi_{n,i}\rightarrow p$\, pour tout  \,$i\in\{1,2\}$\,.
En passant  \`a la limite dans le terme de droite de l'\'{e}galit\'{e}
\eqref{desig:lem:loc:02} on obtient le vecteur unitaire \,$v$.
Par cons\'equent, en revenant \`a l'in\'{e}galit\'{e} 
\eqref{desig:lem:loc} on obtient  la m\^eme in\'{e}galit\'{e}
\eqref{contra:01}
 ce qui produit \`a nouveau, une contradiction.
 Le cas o\`u
\,$\|q_{n}-p_{n}\| \rightarrow 0$\, est similaire au cas \,$\|\lambda_{n}-p_{n}\|\rightarrow 0$\,.
\end{proof}

%%%%%%%%%%%%%%%%%%%%%
\vskip30pt
\section{Pr\'{e}paration pour la preuve du th\'eor\`eme}
\vskip10pt
%%%%%%%%%%%%%%%%%%%%%

Dans la suite de ce travail, le voisinage \,$\mathcal{V}$ de l'identit\'e dans 
\,$\text{Diff}^{1}(\RR)$\, sera donn\'e par
$$\mathcal{V}=\mathcal{U}_{\epsilon} \quad \text{o\`u} \quad 0<\epsilon<\min\{\epsilon_{1}\,,\epsilon_{2}\}$$
et \,$\epsilon_{2}=\epsilon_{2}(\delta,K)$\, est choisi en prenant
\,$\delta=2\epsilon_{1}$\, et \,$K=3$\,. Les voisinages 
\,$\mathcal{U}_{\epsilon_{1}}$\, et \,$\mathcal{U}_{\epsilon_{2}}$\,
ont \'et\'e introduits dans la section \ref{sec:prop:difeo}.

%%%%%%%%%%%%%%%%%%%%%%%%%%%%%%%%%%
%:Observacao: obs:01
\vskip10pt
\begin{remark}\label{obs:01}
Soient  \,$f,h\in \mathcal{V}$\, deux \difeos commutants et soit
\,$x\in\RR$\, tel que \,$f(x)\neq x$.
D'apr\`es le lemme \ref{prop:loc:dif} et le choix du voisinage 
\,$\mathcal{V}$ on a 
\begin{align}\label{}
\big\| \, h(\lambda)-h\big(f^{i}(x)\big)\big\|    &  
\leq 3 \, \big\| \, h\big(f^{i+1}(x) \big)-h\big(f^{i}(x)\big)\big\|
 \notag
\end{align}
pour tout \,$\lambda\in [f^{i}(x)\,,f^{i+1}(x)]$\, et \,$i\in\Z$.
Pourtant, comme \,$f$\, et \,$h$\,  commutent il s'en suit que \,$h(\lambda)$\, est contenu dans la boule ferm\'ee de centre \,$f^{i}\big(h(x)\big)$\, et de rayon
\,$3\,\big\|f^{i+1}\big(h(x)\big)-f^{i}\big(h(x)\big)\big\|$\,.

D'autre part, on sait d'apr\`es \,$(B.1)$\, que \,$f$\, est sans point fixe sur la boule ferm\'ee de centre \,$f^{i}(h(x))$\, et de rayon 
\,$4 \, \big\| \, f^{i+1}\big(h(x) \big)-f^{i}\big(h(x)\big)\big\|$\,. 
Donc, on peut faire une homothopie entre les courbes
$$\big[f^{i}\big(h(x)\big)\,,f^{i+1}\big(h(x)\big)\big] \quad \text{et} \quad 
h\big([f^{i}(x)\,,f^{i+1}(x)]\big)$$ 
sans passer par les points fixes de \,$f$\, et sans bouger les extr\'{e}mit\'{e}s
\,$f^{i}\big(h(x)\big)$\, et \,$f^{i+1}\big(h(x)\big)$\,.

De plus, si \,$(n_{k})_{k\geq1}$\, est une suite croissante dans \,$\Z^{+}$\,
telle que \,$f^{n_{k}}(x)\rightarrow x$\, lorsque \,$k$\, tend vers l'infini et, si on consid\`ere \,$k$\,
suffisamment grand,  alors on peut appliquer aux segments 
\,$[f^{n_{k}}(x)\,,f(x)]$\, le lemme \ref{prop:loc:dif}
ce qui montre que \,$h\big([f^{n_{k}}(x)\,,f(x)]\big)$\, et 
\,$\big[f^{n_{k}}\big(h(x)\big)\,,f\big(h(x)\big)\big]$\,
sont tous les deux dans la boule ouverte de centre \,$h(x)$\, et de rayon
\,$4 \, \big\| \, h(x)-f\big(h(x)\big)\big\| $\, sur laquelle 
\,$f$\, n'a pas de point fixe. Ainsi, on peut, comme pr\'{e}c\'{e}demment, faire une homothopie entre les courbes 
\,$\big[f^{n_{k}}\big(h(x)\big)\,,f\big(h(x)\big)\big]$\,  \text{et}  
\,$h\big([f^{n_{k}}(x)\,,f(x)]\big)$\, 
sans passer par les points fixes de \,$f$\, et sans bouger les extr\'{e}mit\'{e}s
\,$f^{n_{k}}\big(h(x)\big)$\, et \,$f\big(h(x)\big)$\,.

\end{remark}
\vskip5pt
%%%%%%%%%%%%%%%%%%%%%%%%%%%%%%%%%%

Nous allons faire l'usage de cette remarque dans la preuve du lemme
\ref{ptfix:com:fg}.

%%%%%%%%%%%%%%%%%%%%%%%%%%%%%%%%%%
%:Definic‹o: ponto capital
\vskip10pt
\begin{defn}%\label{}
Soit \,$f\in\mathrm{Diff}^{1}(\RR)$\,  et soit
\,${p}\in \RR$\, tels que \,$f({p})\neq {p}$. On dit qu'un point  \,$q\in\mathrm{Fix}(f)$\, est un \,{\it point  capital}\, pour 
\,$\mathcal{O}_{{p}}(f)$\, s'il existe une suite croissante 
\,$(n_{k})_{k\geq1}$\,  dans \,$\Z^{+}$\, telle que:
\begin{itemize}
\item 
$f^{n_{k}}({p})\rightarrow {p}$\, 
lorsque \,$k$\, tend vers l'infini;

\item
$\mathrm{Ind}_{q}(\Gamma^{f}_{{p},n_{k}})$\, est bien d\'efini et non nul pour 
\,$k$\, suffisamment grand.

\end{itemize}

\end{defn}
\vskip5pt
%%%%%%%%%%%%%%%%%%%%%%%%%%%%%%%%%%

Bien s\^ur, si \,$q\in\RR$\, est un point capital
pour \,$\mathcal{O}_{{p}}(f)$\, alors \,$q$\, appartient \`a 
l'enveloppe convexe de \,$\mathcal{O}_{{p}}(f)$\,.

%%%%%%%%%%%%%%%%%%%%%%%%%%%%%%%%%%
%:A limitac‹o da h-—rbita de q: ptfix:com:fg
\vskip10pt
\begin{lemma}\label{ptfix:com:fg}
Soient  \,$f,h\in\mathcal{V}$ deux \difeos commutants et soit
\,$q\in\mathrm{Fix}(f)$\, un point capital pour \,$\mathcal{O}_{p}(f)$\, 
o\`u \,$p\in\RR-\mathrm{Fix}(f)$\,. 
Alors \,$h^{\ell}(q)$\, est un point capital pour 
\,$\mathcal{O}_{h^{\ell}(p)}(f)$\, quel que soit \,$\ell\in\Z$\,
et on a  \,$\mathcal{O}_{q}(h)\subset 
\mathrm{Conv}\big(\mathcal{O}_{p}(f,h)\big)$\,.

\end{lemma}
\vskip5pt
%%%%%%%%%%%%%%%%%%%%%%%%%%%%%%%%%%

\begin{proof}[Preuve]
Comme le point  \,$q$\, est un point capital pour \,$\mathcal{O}_{p}(f)$, il   existe une suite croissante \,$(n_{k})_{k\geq1}$\, 
dans \,$\Z^{+}$\, avec les propri\'et\'es suivantes:
\begin{enumerate}
\item[--] 
$f^{n_{k}}({p})\rightarrow {p}$;

\item[--]
$\text{Ind}_{q}(\Gamma^{f}_{{p},n_{k}})$\, est bien d\'{e}fini et non nul pour 
\,$k$\, suffisamment grand.

\end{enumerate}

On va montrer que les points \,$h(p) \, , \,h(q)$\, et les courbes ferm\'ees
orient\'ees  \,$\Gamma^{f}_{h({p}),n_{k}}$\, ont  les m\^eme  propri\'et\'es qu'\'enonc\'es pr\'ec\'edemment pour   \,$p\,,q $\, et \,$\Gamma^{f}_{{p},n_{k}}$.

Comme \,$f$\, et \,$g$\, commutent on a 
\,$h(q)\in\text{Fix}(f)$\, et  on aura encore:
 $$h\big(\mathcal{O}_{p}(f)\big)=\mathcal{O}_{h(p)}(f) \quad \text{et} \quad
f^{n_{k}}\big(h({p})\big)= h\big(f^{n_{k}}({p})\big)\rightarrow h({p})\,.$$

En outre, l'indice de \,$h(q)$\, par rapport \`a la courbe 
 \,$h\big(\Gamma^{f}_{{p},n_{k}}\big)$\, est bien d\'{e}fini et est non nul. D'autre part, il suit de la Remarque  \ref{obs:01}  
qu'on peut faire une homothopie entre 
\,$h(\Gamma^{f}_{{p},n_{k}})$\, et 
\,$\Gamma^{f}_{h({p}),n_{k}}$\,, qui fixe les sommets de 
\,$\Gamma^{f}_{h({p}),n_{k}}$\, et ne passe pas par les points fixes de
\,$f$. On en d\'eduit alors que  
\,$\text{Ind}_{h(q)}\big(\Gamma^{f}_{h({p}),n_{k}}\big)$\,
est bien d\'{e}fini et
$$\text{Ind}_{h(q)}\big(\Gamma^{f}_{h({p}),n_{k}}\big)
=\text{Ind}_{h(q)}\big(h(\Gamma_{{p},n_{k}}^{f})\big)\neq 0$$
pourvu que \,$k$\, soit suffisamment grand. Ainsi, on vient de montrer que 
\,$h(q)$\, est un point capital pour  \,$\mathcal{O}_{h({p})}(f)$\,.
Donc  \,$h(q)$\, appartient \`a l'enveloppe convexe de 
 \,$\mathcal{O}_{h({p})}(f)$\, et, par cons\'equent, on aura 
 \,$h(q)\in\text{Conv}\big(\mathcal{O}_{p}(f,h)\big)$\,.
 
 \vskip10pt

\setlength{\unitlength}{1mm}%
\noindent
\hfil
%\fbox{%
\begin{picture}(58,45)(-10,-5)%

\put(0,0){\includegraphics[width=4cm]{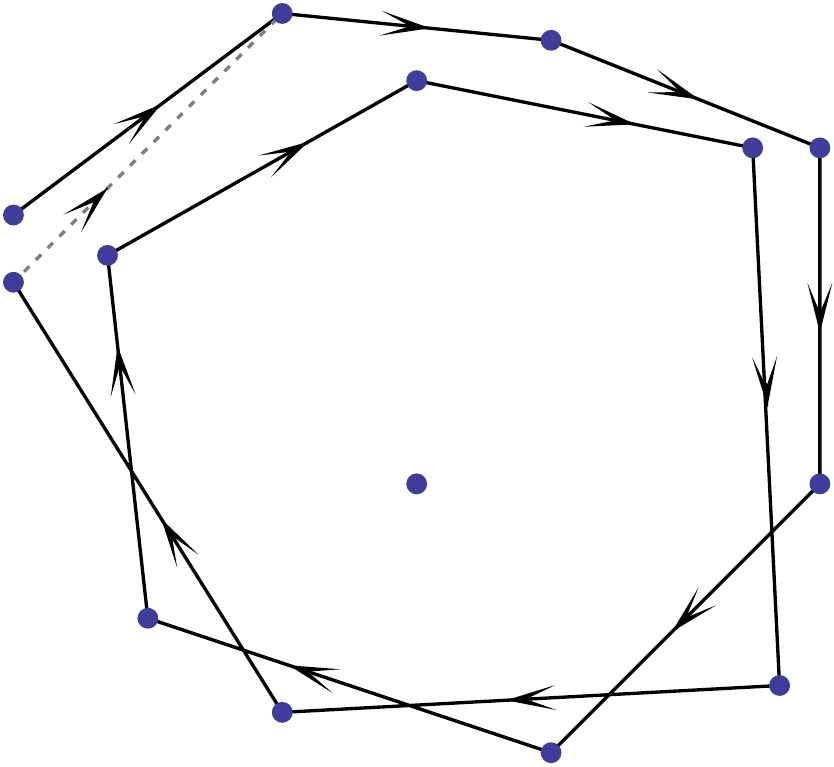}}

\put(-2,27){\tiny${p}$}
\put(9,38){\tiny$f({p})$}
\put(26,36.5){\tiny$f^{2}({p})$}
\put(40,29){\tiny$f^{3}({p})$}
\put(40,13){\tiny$f^{4}({p})$}
\put(25,-3){\tiny$f^{5}({p})$}
\put(-2,6){\tiny$f^{6}({p})$}
\put(-10,22){\tiny$f^{n_{k}}({p})$}
\put(18,15){\tiny$q$}

\put(40,3){\scriptsize$\Gamma^{f}_{{p},n_{k}}$}
\put(38,5){$\leftarrow$}

\end{picture}%
%}
\hfil
\setlength{\unitlength}{1mm}%
\noindent
\hfil
%\fbox{%
\begin{picture}(66,45)(-14,-5)%

\put(0,0){\includegraphics[width=4cm]{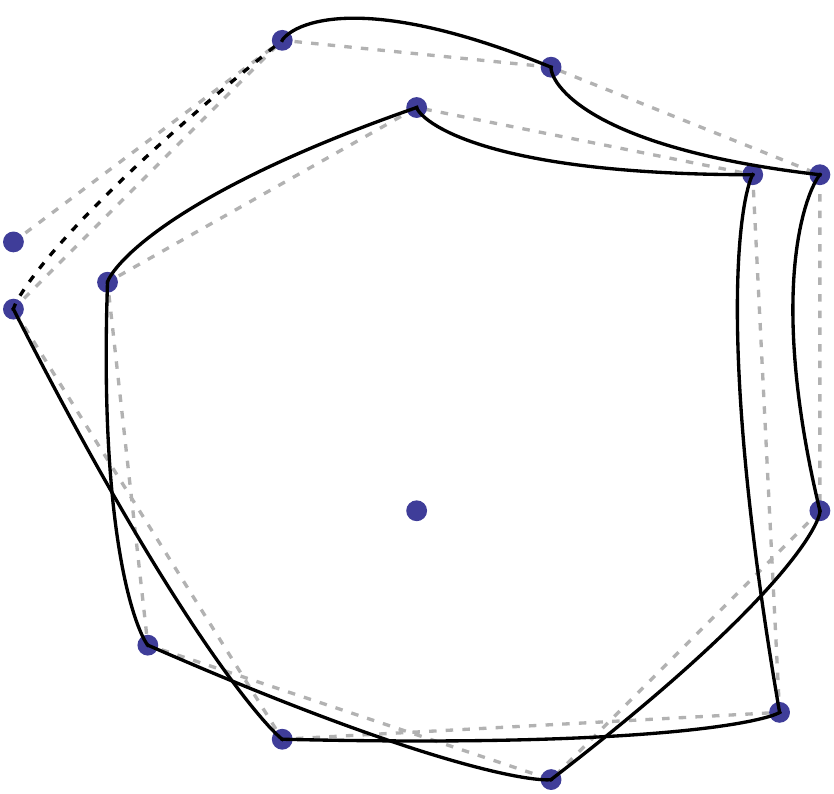}}

\put(-5,28){\tiny$h({p})$}
\put(2,36){\tiny$f(h({p}))$}
\put(26,36.5){\tiny$f^{2}(h({p}))$}
\put(40,29){\tiny$f^{3}(h({p}))$}
\put(40,13){\tiny$f^{4}(h({p}))$}
\put(25,-3){\tiny$f^{5}(h({p}))$}
\put(-6,6){\tiny$f^{6}(h({p}))$}
\put(-14,22){\tiny$f^{n_{k}}(h({p}))$}
\put(18,11){\tiny$h(q)$}

\put(20.5,23){\scriptsize$h(\Gamma^{f}_{{p},n_{k}})\!\rightarrow$}
\put(39,23){\scriptsize$\leftarrow\Gamma^{f}_{h({p}),n_{k}}$}

\end{picture}%
%}

En r\'{e}p\'{e}tant ce processus successivement on conclut que:
\begin{enumerate}
\item[--] 
pour chaque \,$\ell\in\Z^{+}$, l'indice 
\,$\mathrm{Ind}_{h^{\ell}(q)}\big(\Gamma^{f}_{h^{\ell}({p}),n_{k}}\big)$\, 
est bien d\'{e}fini et non nul lorsque \,$k$\, est suffisamment grand;

\item[--]
$h^{\ell}(q)\in \text{Conv}\big({\mathcal{O}_{p}(f,h)}\big)$\, pour tout entier
\,$\ell\geq0$.

\end{enumerate}

La preuve pour le cas \,$\ell<0$\, est similaire, il suffit de diminuer le voisinage \,$\mathcal{V}$ pour qu'on puisse appliquer les raisonnements \`a \,$h^{-1}$. Ceci cl\^ot la preuve du lemme. 
\end{proof}

%%%%%%%%%%%%%%%%%%%%%%%%%%%%%%%%%%
%:O Teorema no caso n=1
\vskip10pt
\begin{lemma}\label{prep02:exis:ptf:comu}
Soient  \,$g_{1},\ldots,g_{m}\,,f\in \mathcal{V}$\, des \difeos commutants et soit
\,$p\in\mathrm{Fix}(g_{1},\ldots,g_{m})$\, un point dont la \,$f$-orbite 
est born\'ee. Alors, il existe un point 
\,$q\in\mathrm{Fix}(g_{1},\ldots,g_{m}\,;f)$\, 
satisfaisant une des conditions$\,:$
\begin{enumerate}
\item[$(i)$]
$q\in \overline{\mathcal{O}_{p}(f)}\,;$

\item[$(ii)$]
$q\in \mathrm{Conv}\big(\mathcal{O}_{p}(f)\big)$\, et dans ce cas \,$q$\,
est un point capital pour \,$\mathcal{O}_{\tilde{p}}(f)$\, o\`u 
\,$\tilde{p}\in \overline{\mathcal{O}_{p}(f)}-\mathrm{Fix}(f)$.
\end{enumerate}

\end{lemma}
\vskip5pt
%%%%%%%%%%%%%%%%%%%%%%%%%%%%%%%%%%

\begin{proof}[Preuve]
Admettons  qu'il n'y ait pas de points fixes de \,$f$\, dans 
\,$\overline{\mathcal{O}_{p}(f)}$\,. Ainsi,
\,$\overline{\mathcal{O}_{p}(f)}\subset\text{Fix}(g_{1},\ldots,g_{m})$\,
contient un compact minimal \,$\Lambda$\, sans point fixe pour \,$f$.
Dans ce cas, \'etant donn\'e 
\,$\tilde{p}\in\Lambda \subset \overline{\mathcal{O}_{p}(f)}$\,,
on a  \,$\overline{\{f^{i}(p) \ ; \ i\geq0\}}=\Lambda$\,.
Donc, il existe une suite croissante \,$(n_{k})_{k\geq1}$\,
dans \,$\Z^{+}$ telle que 
\,$f^{n_{k}}(\tilde{p})\rightarrow \tilde{p}$\,.

Maintenant, passons  aux courbes ferm\'ees et orient\'ees 
\,$\Gamma^{f}_{\tilde{p},n_{k}}$\,.
D'apr\`es le lemme \ref{lema:fundam:ind}, il existe
\,$q_{k}\in \text{Fix}(g_{1},\ldots,g_{m}\,;f)$\, 
tel que \,$\text{Ind}_{q_{k}}(\Gamma^{f}_{\tilde{p},n_{k}})$\, est bien d\'{e}fini et non nul pourvu que \,$k$\, soit suffisamment grand. Et dans ces conditions 
\,$q_{k}$\, appartient \`a \,$\mathrm{Conv}\big(\mathcal{O}_{\tilde{p}}(f)\big)$\,. 

Comme la \,$f$-orbite de \,$p$\, est born\'ee et  
\,$\tilde{p}\in \overline{\mathcal{O}_{p}(f)}$\,, on en d\'eduit que
la \,$f$-orbite de \,$\tilde{p}$\, est born\'ee et donc,
\,$\mathrm{Conv}\big(\mathcal{O}_{\tilde{p}}(f)\big)$\, est aussi
born\'e. Dans ce cas, quitte \`a extraire une sous suite de 
\,$(q_{k})_{k\geq1}$\,,  on peut admettre sans perte de g\'{e}n\'{e}ralit\'{e}, que \,$q_{k}\rightarrow q$\, avec \,$q$\, dans \,$\text{Fix}(g_{1},\ldots,g_{m}\,;f)$.

Comme \,$\overline{\mathcal{O}_{p}(f)}\cap \text{Fix}(f)=\emptyset$\,, il existe une boule 
\,$B(q\,;\delta)$\, telle que 
\,$B(q\,;\delta) \cap \overline{\mathcal{O}_{p}(f)}=\emptyset$\,.

Supposons qu'on puisse trouver des segments \,$[p',f(p')]$\,   arbitrairement proches de 
\,$q$\, avec \,$p'\in\overline{\mathcal{O}_{p}(f)}$\,. Dans ce cas, 
 la compacit\'e de \,$\overline{\mathcal{O}_{p}(f)}$\, montre
qu'il existe \,$q'\in \overline{\mathcal{O}_{p}(f)}-\text{Fix}(f)$\, tel que 
\,$q\in[q',f(q')]$\,; ceci donne 
 une contradiction car, par \,$(B.1)$\,, il n'existe pas de points fixes de \,$f$\, dans le segment \,$[q',f(q')]$. De plus, en prenant \,$k$\, suffisamment grand, on d\'eduit que les segments \,$[f^{n_{k}}(\tilde{p}),f(\tilde{p})]$\, ne peuvent pas passer arbitrairement proches 
de \,$q$\,  car \,$f^{n_{k}}(\tilde{p})\rightarrow \tilde{p}$\, et 
\,$[\tilde{p},f(\tilde{p})]$\, ne contient pas le point \,$q$. 

\vskip10pt

\setlength{\unitlength}{1mm}%
\noindent
\hfil
%\fbox{%
\begin{picture}(58,45)(-10,-5)%

\put(0,0){\includegraphics[width=4cm]{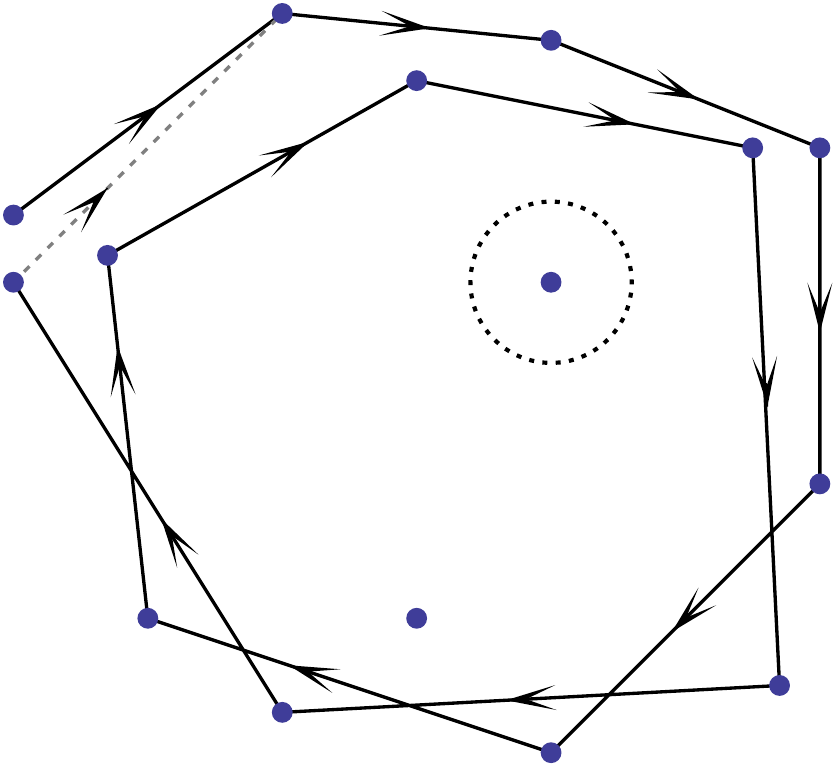}}

\put(-2,27){\tiny$\tilde{p}$}
\put(9,38){\tiny$f(\tilde{p})$}
\put(26,36.5){\tiny$f^{2}(\tilde{p})$}
\put(40,29){\tiny$f^{3}(\tilde{p})$}
\put(40,13){\tiny$f^{4}(\tilde{p})$}
\put(25,-3){\tiny$f^{5}(\tilde{p})$}
\put(-2,6){\tiny$f^{6}(\tilde{p})$}
\put(-10,22){\tiny$f^{n_{k}}(\tilde{p})$}
\put(19,5){\tiny$q_{k}$}
\put(25,21){\tiny$q$}
\put(13,23){\tiny$B(q\,;\delta)$}
\put(40,3){\scriptsize$\Gamma^{f}_{\tilde{p},n_{k}}$}
\put(38,5){$\leftarrow$}

\end{picture}%
%}

Alors, quitte \`a passer
\`a une valeur plus petite de \,$\delta$\,, on peut supposer que la boule 
\,$B(q\,;{\delta})$\, isole le point \,$q$\, de toutes les courbes 
\,$\Gamma^{f}_{\tilde{p},n_{k}}$\, si  \,$k$\, est assez grand.
Et dans cette condition,   l'indice de \,$q$\, par rapport \`a
la courbe \,$\Gamma^{f}_{\tilde{p},n_{k}}$\, est bien d\'{e}fini et on obtient:
\begin{enumerate}
\item[--] 
$q\in\text{Fix}(g_{1},\ldots,g_{m}\,;f)$\,;

\item[--]
$\text{Ind}_{q}(\Gamma^{f}_{\tilde{p},n_{k}})=
\text{Ind}_{q_{k}}(\Gamma^{f}_{\tilde{p},n_{k}})\neq 0 $\,, 
\end{enumerate}
ce qui d\'emontre l'existence d'un point capital
\,$q\in\text{Fix}(g_{1},\ldots,g_{m}\,;f)$\, associ\'e \`a
\,$\mathcal{O}_{\tilde{p}}(f)$\, o\`u 
\,$\tilde{p}\in\overline{\mathcal{O}_{p}(f)}$.

Pour conclure la preuve, il suffit de se rappeler que, en prenant
\,$n_{k}$\, grand et 
\,$\tilde{p}'\in\mathcal{O}_{p}(f)$\,  proche de \,$\tilde{p}$\,,
 l'indice \,$\text{Ind}_{q}(\Gamma_{\tilde{p}'\!,n_{k}}^{f})$\, 
est bien d\'{e}fini et il co\"{\i}ncide avec 
\,$\text{Ind}_{q}(\Gamma_{\tilde{p},n_{k}}^{f})$\,. Alors, 
\,$q\in \text{Conv}\big(\mathcal{O}_{p}(f)\big)$\, et le lemme est
conclut.
\end{proof}

%%%%%%%%%%%%%%%%%%%%%
\vskip30pt
\section{Preuve du th\'eor\`eme \ref{teor:central}}
\vskip10pt
%%%%%%%%%%%%%%%%%%%%%

Pour d\'emontrer le th\'eor\`eme \ref{teor:central}, on commence par le cas o\`u il y a un nombre fini de g\'{e}n\'{e}rateurs.

%%%%%%%%%%%%%%%%%%%%%%%%%%%%%%%%%
\vskip10pt
\begin{thm}\label{teo:cas:fini}
Soient   \,$g_{1},\ldots,g_{m}\,,f_{1},\ldots,f_{n}\in\mathcal{V}$
des \difeos commutants et supposons qu'il existe un point 
\,$p\in\mathrm{Fix}(g_{1},\ldots,g_{m})$\, dont la 
\,$(f_{1},\ldots,f_{n})$-orbite est born\'ee. Alors, il existe un point  fixe 
\,$q$\, commun \`a \,$g_{1},\ldots,g_{m}\,,f_{1},\ldots,f_{n}$\, 
satisfaisant une des conditions$\,:$
\begin{enumerate}
\item [$(i)$]
$q\in\overline{\mathcal{O}_{p}(f_{1},\ldots,f_{n})} \,;$

\item[$(ii)$]
$q\in\mathrm{Conv}\Big({\mathcal{O}_{p}(f_{1},\ldots,f_{n})}\Big)$\, 
et dans ce cas \,$q$\, est un point capital pour
\,$\mathcal{O}_{\tilde{p}}(f_{j})$\, o\`u 
\,$j\in\{1,\ldots,n\}$\, et 
\,$\mathcal{O}_{\tilde{p}}(f_{j})\subset \mathrm{Conv}\Big(\,\overline{\mathcal{O}_{p}(f_{1},\ldots,f_{n})}\,\Big)$\,.

\end{enumerate}

\end{thm}
\vskip5pt
%%%%%%%%%%%%%%%%%%%%%%%%%%%%%%%%%

\begin{proof}[Preuve]
On fera une preuve par r\'ecurrence sur \,$n$. Le lemme \ref{prep02:exis:ptf:comu} 
donne le r\'esultat dans le cas \,$n=1$\, et \,$m\in\Z^{+}$.
Admettons alors que le r\'{e}sultat est vrai pour  \,$n\geq1$\, quelconque et
d\'emontrons le dans le cas \,$n+1$.

Pour cela, soient  \,$g_{1},\ldots,g_{m}\,,f_{1},\ldots,f_{n}\,,f_{n+1}\in\mathcal{V}$ des \difeos commutants et soit \,$p\in \text{Fix}(g_{1},\ldots,g_{m})$\,
dont la \,$(f_{1},\ldots,f_{n}\,,f_{n+1})$-orbite est born\'ee.

Le lemme \ref{prep02:exis:ptf:comu} appliqu\'e au point \,$p$\, et aux 
\difeos \,$g_{1},\ldots,g_{m}\,,f_{1}$\,  donne un point fixe \,$p'$  commun aux \difeos \,$g_{1},\ldots,g_{m}\,,f_{1}$\, et
satisfaisant l'une des conditions$\,:$

\begin{enumerate}
\item[\it Cas I\,$:$]
{\it $p'\in \overline{\mathcal{O}_{p}(f_{1})} \,;$}

\item[\it Cas II\,$:$]
{\it $p' \in\text{\rm Conv}\big(\mathcal{O}_{p}(f_{1})\big)$\, et dans ce cas
\,$p'$ est point capital pour 
\,$\mathcal{O}_{\tilde{p}}(f_{1})$\, o\`u \,$\tilde{p}\in
\overline{\mathcal{O}_{p}(f_{1})}$}\,. 

\end{enumerate}

\'Etudions le \emph{Cas I}. Par hypoth\`ese,
 la \,$(f_{2},\ldots,f_{n}\,,f_{n+1})$-orbite de \,$p'$ est contenue dans
\,$\overline{\mathcal{O}_{p}(f_{1},\ldots,f_{n}\,,f_{n+1})}$\, et donc, est born\'ee. L'hypoth\`ese d'induction appliqu\'ee au point 
\,$p'\in\text{Fix}(g_{1},\ldots,g_{m}\,;f_{1})$\,  et aux \difeos
$$g_{1},\ldots,g_{m}\,,f_{1} \quad \text{e} \quad f_{2},\ldots,f_{n}\,,f_{n+1}$$
nous garantit qu'il existe un point fixe
$$q\in \text{Fix}(g_{1},\ldots,g_{m}\,,f_{1}\,;f_{2},\ldots,f_{n}\,,f_{n+1})$$ 
satisfaisant l'une des conditions$\,:$

\begin{enumerate}
\item[\it Cas I-A\,$:$]
$q\in \overline{\mathcal{O}_{p'}(f_{2},\ldots,f_{n}\,,f_{n+1})} \,;$

\item[\it Cas I-B\,$:$]
{\it $q\in \text{\rm Conv}\Big({\mathcal{O}_{p'}(f_{2},\ldots,f_{n}\,,f_{n+1})}\Big)$\,  et dans ce cas \,$q$\, est un point capital pour
\,$\mathcal{O}_{\tilde{p}'}(f_{j})$\, o\`u
\,$j\in\{2,\ldots,n\,,n+1\}$\, et 
$$\mathcal{O}_{\tilde{p}'}(f_{j})\subset \text{\rm Conv}\Big(\,\overline{\mathcal{O}_{p'}(f_{2},\ldots,f_{n}\,,f_{n+1})}\,\Big) \,.$$}
\end{enumerate}

\vskip5pt

Si on se place sur le \,{\it Cas I-A}\,,
alors \,$q\in \overline{\mathcal{O}_{p'}(f_{1},f_{2},\ldots,f_{n},f_{n+1})}$\,
car \,$p'\in\text{Fix}(f_{1})$\, et, comme  
\,$p'\in\overline{\mathcal{O}_{p}(f_{1})}$\,, on doit avoir
$$q\in \overline{\mathcal{O}_{p}(f_{1},f_{2},\ldots,f_{n},f_{n+1})}\,.$$

\vskip5pt

Supposons maintenant que le \,{\it Cas I-A}\, n'arrive pas. Ainsi,
$$q\in\text{Fix}(g_{1},\ldots,g_{m}\,,f_{1}\,;f_{2},\ldots,f_{n},f_{n+1}) $$
et est contenu dans l'enveloppe convexe de 
\,${\mathcal{O}_{p'}(f_{2},\ldots,f_{n},f_{n+1})}$. Donc,
$$q\in\text{Conv}\Big(\,\overline{\mathcal{O}_{p}(f_{1},f_{2},\ldots,f_{n},f_{n+1})}\,\Big)\,.$$
En plus, le point \,$q$ est un point capital pour 
\,$\mathcal{O}_{\tilde{p}'}(f_{j})$\,
o\`u \,$j\in\{2,\ldots,n\,,n+1\}$\, et 
$$\mathcal{O}_{\tilde{p}'}(f_{j})\subset\text{Conv}\Big(\,\overline{\mathcal{O}_{p'}(f_{2},\ldots,f_{n},f_{n+1})}\,\Big)$$
et donc,
$$\mathcal{O}_{\tilde{p}'}(f_{j})\subset\text{Conv}\Big(\,\overline{\mathcal{O}_{p}(f_{1},f_{2},\ldots,f_{n},f_{n+1})}\,\Big)$$
car \,$p'\in\overline{\mathcal{O}_{p}(f_{1})}$\,. 

Il reste \`a prouver 
\,$q\in\text{Conv}\Big({\mathcal{O}_{p}(f_{1},f_{2},\ldots
,f_{n},f_{n+1})}\Big)\,.$
Rappelons-nous que 
$$\text{Conv}\Big(\,\overline{\mathcal{O}_{p}(f_{1},f_{2},\ldots,f_{n},f_{n+1})}\,\Big)=\overline{\text{Conv}\Big({\mathcal{O}_{p}(f_{1},f_{2},\ldots,f_{n},f_{n+1})}\,\Big)}\,.$$
Ainsi,  on peut approcher les sommets de la courbe \,$\Gamma^{f_{j}}_{\tilde{p}',n_{k}}$\, pour des points qui sont dans
\,$\text{Conv}\Big({\mathcal{O}_{p}(f_{1},f_{2},\ldots,f_{n},f_{n+1})}\,\Big)$\, 
et construire une courbe 
$$\gamma=\Gamma(a_{1}, \ldots,a_{\ell})$$ 
qui ne passe pas par le point \,$q$\, et telle que
\,$\text{Ind}_{q}(\gamma)=\text{Inf}_{q}\big(\Gamma^{f_{j}}_{\tilde{p}',n_{k}}\big)\neq0$\,. Cela d\'emontre que le point \,$q$\, est contenu dans 
\,$\text{Conv}\Big({\mathcal{O}_{p}(f_{1},f_{2},\ldots,f_{n},f_{n+1})}\,\Big)$\, et la preuve est finie pour le \,{\it Cas I}\,.

\vskip5pt

Supposons maintenant que le \,{\it Cas I}\, n'arrive pas. Alors, le point 
$$p'\in \text{Fix}(g_{1},\ldots,g_{m}\,;f_{1})$$ 
est un point capital pour \,$\mathcal{O}_{\tilde{p}}(f_{1})$\, o\`u
\,$\tilde{p}\in\overline{\mathcal{O}_{p}(f_{1})}$.
D'autre part, le lemme \ref{ptfix:com:fg} appliqu\'e successivement aux 
\difeos
$$f_{1},h \quad \text{avec} \quad h=f_{2},\ldots,f_{n}\,,f_{n+1}$$
nous garantit que la \,$(f_{2},\ldots,f_{n}\,,f_{n+1})$-orbite
du point \,$p'$\, est contenue dans l'enveloppe convexe de 
\,$\mathcal{O}_{\tilde{p}}(f_{1},\ldots,f_{n}\,,f_{n+1})$\, et, par cons\'equent,
est born\'ee puisque
$$\mathcal{O}_{\tilde{p}}(f_{1},\ldots,f_{n}\,,f_{n+1}) \subset 
\overline{\mathcal{O}_{{p}}(f_{1},\ldots,f_{n}\,,f_{n+1})}$$
et \,$\tilde{p}\in\overline{\mathcal{O}_{p}(f_{1})}$\,.
On peut alors appliquer l'hypoth\`ese d'induction au point 
\,$p'\in \text{Fix}(g_{1},\ldots,g_{m}\,;f_{1})$\, et aux \difeos
$$g_{1},\ldots,g_{m}\,,f_{1} \quad \text{et} \quad f_{2},\ldots,f_{n}\,,f_{n+1}$$
pour conclure qu'il existe un point fixe \,$q$\, commun aux \difeos
$$g_{1},\ldots,g_{m}\,,f_{1},f_{2},\ldots,f_{n}\,,f_{n+1}$$
satisfaisant l'une des conditions$\,:$

\begin{enumerate}
\item[\it Cas II-A\,$:$]
$q\in \overline{\mathcal{O}_{p'}(f_{2},\ldots,f_{n}\,,f_{n+1})} \,;$

\item[\it Cas II-B\,$:$]
$q\in \text{Conv}\Big({\mathcal{O}_{p'}(f_{2},\ldots,f_{n}\,,f_{n+1})}\Big)$\,
{\it  et dans ce cas \,$q$\, est un point ca\-pi\-tal pour  \,$\mathcal{O}_{\tilde{p}'}(f_{j})$\, o\`u
\,$j\in\{2,\ldots,n\,,n+1\}$\, et}\, 
$$\mathcal{O}_{\tilde{p}'}(f_{j})\in \text{Conv}\Big(\,\overline{\mathcal{O}_{p'}(f_{2},\ldots,f_{n}\,,f_{n+1})}\,\Big).$$
\end{enumerate}

Supposons que le \,{\it Cas II-A}\, arrive. Si le point \,$q$\, est dans  
\,$\overline{\mathcal{O}_{p}(f_{1},\ldots,f_{n}\,,f_{n+1})}$\,,
alors la preuve est finie.
Par contre, si 
\,$q\notin\overline{\mathcal{O}_{p}(f_{1},\ldots,f_{n}\,,f_{n+1})}$\, on en  d\'eduit qu'il existe une boule  \,$B(q\,;\delta)$\, telle que 
$$B(q\,;\delta)\cap\overline{\mathcal{O}_{p}(f_{1},\ldots,f_{n}\,,f_{n+1})}=\emptyset \,.$$

Maintenant, admettons qu'il existe des segments \,$[\tilde{q}\,,f_{1}(\tilde{q})]$\, passant arbitrairement proches du point \,$q$\, avec
\,$\tilde{q}\in \overline{\mathcal{O}_{p}(f_{1},\ldots,f_{n}\,,f_{n+1})}$\,.
Il suit de la compacit\'e de 
\,$\overline{\mathcal{O}_{p}(f_{1},\ldots,f_{n}\,,f_{n+1})}$\,
qu'il existe un point 
\,$\tilde{q}'\in\overline{\mathcal{O}_{p}(f_{1},\ldots,f_{n}\,,f_{n+1})}$\,
tel que \,$q\in [\tilde{q}',f_{1}(\tilde{q}')]$\, ce qui est  un absurde
car il n'existe pas de points fixes de 
\,$f_{1}$\, sur \,$[\tilde{q}',f_{1}(\tilde{q}')]$\,. 
Ainsi, en prenant \,$\delta$\, plus petit, si n\'{e}cessaire, on peut admettre
que \,$B(q\,;\delta)$\, isole \,$q$\, de l'ensemble form\'e par toutes les courbes \,$\Gamma_{\tilde{q}}^{f_{1}}$\, o\`u
\,$\tilde{q}\in \overline{\mathcal{O}_{p}(f_{1},\ldots,f_{n}\,,f_{n+1})}$\,.

En plus, comme 
\,$q\in \overline{\mathcal{O}_{p'}(f_{2},\ldots,f_{n}\,,f_{n+1})}$\, il s'en suit
qu'il existe des suites d'entiers
\,$(\xi_{2,\ell})_{\ell\geq1}\,, \ldots , (\xi_{n,\ell})_{\ell\geq1}\,,
(\xi_{n+1,\ell})_{\ell\geq1}$\,
tels que 
$$f_{2}^{\xi_{2,\ell}}\circ \cdots \circ f_{n}^{\xi_{n,\ell}}\circ
f_{n+1}^{\xi_{n+1,\ell}}(p') \rightarrow q .$$
Dans la suite, on va utiliser la notation
$$h_{\ell}:=f_{2}^{\xi_{2,\ell}}\circ \cdots \circ f_{n}^{\xi_{n,\ell}}\circ
f_{n+1}^{\xi_{n+1,\ell}}\,.$$

\vskip10pt

\setlength{\unitlength}{1mm}%
\noindent
\hfil
%\fbox{%
\begin{picture}(58,45)(-10,-5)%

\put(0,0){\includegraphics[width=4cm]{PontoEssencial-01}}

\put(-2,27){\tiny$\tilde{p}$}
\put(9,38){\tiny$f_{1}(\tilde{p})$}
\put(26,36.5){\tiny$f^{2}_{1}(\tilde{p})$}
\put(40,29){\tiny$f^{3}_{1}(\tilde{p})$}
\put(40,13){\tiny$f^{4}_{1}(\tilde{p})$}
\put(20,-3){\tiny$f^{5}_{1}(\tilde{p})$}
\put(-2,6){\tiny$f^{6}_{1}(\tilde{p})$}
\put(-10,22){\tiny$f^{n_{k}}_{1}(\tilde{p})$}
\put(19,11){\tiny$p'$}

\put(35,1){\scriptsize$\nwarrow\Gamma^{f}_{\tilde{p},n_{k}}$}

\end{picture}%
%}
\hfil
\setlength{\unitlength}{1mm}%
\noindent
\hfil
%\fbox{%
\begin{picture}(67,45)(-15,-5)%

\begin{rotate}{-30}
\put(-15,5){\includegraphics[width=4cm]{PontoEssencial-03}}
\end{rotate}

\put(-2,36){\tiny$h_{\ell}(\tilde{p})$}
\put(16,38){\tiny$f_{1}(h_{\ell}(\tilde{p}))$}
\put(31,29){\tiny$f^{2}_{1}(h_{\ell}(\tilde{p}))$}
\put(39,17){\tiny$f^{3}_{1}(h_{\ell}(\tilde{p}))$}
\put(31,3){\tiny$f^{4}_{1}(h_{\ell}(\tilde{p}))$}
\put(24,-5){\tiny$f^{5}_{1}(h_{\ell}(\tilde{p}))$}
\put(-1,-3){\tiny$f^{6}_{1}(h_{\ell}(\tilde{p}))$}
\put(-15,30){\tiny$f^{n_{k}}_{1}(h_{\ell}(\tilde{p}))$}
\put(8,6){\tiny$h_{\ell}(p')$}
\put(23,16){\tiny$q$}
\put(15,23){\tiny$B(q\,;\delta)$}

\put(-12,23){\scriptsize$\Gamma^{f}_{h_{\ell}(\tilde{p}),n_{k}}\!\!
\rightarrow$}

\end{picture}%
%}

\vskip10pt

Rappelons-nous  que d'apr\`es lemme \ref{ptfix:com:fg}, on sait que 
\,$h_{\ell}(p')$\, est un point capital pour 
\,$\mathcal{O}_{h_{\ell}(\tilde{p})}(f_{1})$\, o\`u 
\,$\tilde{p}\in \overline{\mathcal{O}_{p}(f_{1})}$\,.

Admettons maintenant qu'on a des segments du type
\,$\big[f_{1}^{n_{k}}\big(h_{\ell}(\tilde{p})\big),f_{1}\big(h_{\ell}(\tilde{p})\big)\big]$\,
passant arbitrairement proches de \,$q$\, pour \,$k$\, et \,$\ell$\, grands.
Quitte \`a passer \`a une sous-suite de
\,$(h_{\ell})_{\ell\geq1}$\,, on peut admettre que 
\,$h_{\ell}(\tilde{p})\rightarrow \tilde{p}'$
car \,$\mathcal{O}_{p}(f_{1},\ldots,f_{n}\,,f_{n+1})$\,
est born\'ee et \,$\tilde{p}\in\overline{\mathcal{O}_{p}(f_{1})}$.
D'autre part, pour chaque \,$h_{\ell}(\tilde{p})$\, fix\'e, on sait que 
\,$f_{1}^{n_{k}}\big(h_{\ell}(\tilde{p})\big)\rightarrow h_{\ell}(\tilde{p})$\,. 
Ainsi, on peut construire une sous-suite 
\,$\Big(f_{1}^{n_{k_{s}}}\big(h_{\ell_{s}}(\tilde{p})\big)\Big)_{s\geq1}$\,
de \,$\Big(f_{1}^{n_{k}}\big(h_{\ell}(\tilde{p})\big)\Big)_{\ell\geq1}$\, 
telle que 
\,$f_{1}^{n_{k_{s}}}\big(h_{\ell_{s}}(\tilde{p})\big)\rightarrow \tilde{p}'$\,
et on a en plus  les segments 
\,$\big[f_{1}^{n_{k_{s}}}\big(h_{\ell_{s}}(\tilde{p})\big),f_{1}\big(h_{\ell_{s}}(\tilde{p})\big)\big]$\, qui passent arbitrairement proches de \,$q$\, pour 
\,$s$\, grand.

Dans ce cas, le segment
\,$\Big[f_{1}^{n_{k_{s}}}\big(h_{\ell_{s}}(\tilde{p})\big),f_{1}\big(h_{\ell_{s}}(\tilde{p})\big)\Big]$\,
s'approche arbitrairement du segment \,$[\tilde{p}',f_{1}(\tilde{p}')]$\, et pourtant \,$[\tilde{p}',f_{1}(\tilde{p}')]$\, contient le point \,$q$\, ce qui
est un absurde.
Alors, quitte \`a diminuer \,$\delta>0$\, on peut admettre que la boule
\,$B(q\,;\delta)$\, isole  le point \,$q$\, de toutes les courbes 
\,$\Gamma^{f_{1}}_{h_{\ell}(\tilde{p}),n_{k}}$\, lorsque 
\,$k$\, et \,$\ell$\, sont grands.

Maintenant, en prenant \,$\ell$\, suffisamment grand, on  d\'eduit que 
le point \,$q$\, est un point capital pour 
\,$\mathcal{O}_{h_{\ell}(\tilde{p})}(f_{1})$\,
o\`u  \,$\tilde{p}\in\overline{\mathcal{O}_{p}(f_{1})}$\,. 
Donc, le point  \,$q$\, est contenu dans l'enveloppe convexe de 
\,$\mathcal{O}_{p}(f_{1},\ldots,f_{n}\,,f_{n+1})$\, parce que
\,$h_{\ell}(\tilde{p})\in\overline{\mathcal{O}_{h_{\ell}({p})}(f_{1})}$\,.
En outre, on  d\'eduit que 
\,$\mathcal{O}_{h_{\ell}(\tilde{p})}(f_{1}) \subset 
\overline{\mathcal{O}_{p}(f_{1},\ldots,f_{n}\,,f_{n+1})}$\, et le th\'eor\`eme 
vient d'\^etre prouv\'e lorsque le \,{\it Cas II-A}\, se produit.

\vskip5pt

Supposons maintenant que le \,{\it Cas II-A}\, n'arrive pas. Ainsi, le point 
\,$q$ commun aux \difeos 
$$g_{1},\ldots,g_{m},f_{1},\ldots,f_{n},f_{n+1}$$ appartient \`a
\,$\text{Conv}\Big({\mathcal{O}_{p'}(f_{2},\ldots,f_{n},f_{n+1})}\Big)$\, et 
est un point capital pour \,$\mathcal{O}_{\tilde{p}'}(f_{j})$\, avec \,$j\in\{2,\ldots,f_{n}\,,f_{n+1}\}$\, et 
\,$\mathcal{O}_{\tilde{p}'}(f_{j})\subset \text{Conv}\Big(\overline{\mathcal{O}_{p'}(f_{2},\ldots,f_{n},f_{n+1})}\,\Big)$\,  
o\`u $$p'\in\text{Conv}\big(\mathcal{O}_{p}(f_{1})\big)
\cap \text{Fix}(g_{1},\ldots,g_{m}\,;f_{1}) \,.$$ 

Comme la \,$(f_{2},\ldots,f_{n},f_{n+1})$-orbite de \,$p'$ est contenue dans
$$\text{Conv}\big(\mathcal{O}_{p}(f_{1},\ldots,f_{n},f_{n+1})\big)$$
on  d\'eduit imm\'ediatement que \,$q\in \text{Conv}\Big({\mathcal{O}_{p}(f_{1},\ldots,f_{n},f_{n+1})}\Big)$\,. De plus, il s'en suit pour les m\^emes raisons que
\,$\mathcal{O}_{\tilde{p}'}(f_{j}) 
\subset \text{Conv}\Big(\,\overline{\mathcal{O}_{p}(f_{1},\ldots,f_{n},f_{n+1})}\,\Big)$\, ce qui d\'emontre que le point fixe \,$q$\, est un point capital pour
\,$\mathcal{O}_{\tilde{p}'}(f_{j})$\, avec 
$$\mathcal{O}_{\tilde{p}'}(f_{j})\subset \text{Conv}\Big(\,\overline{\mathcal{O}_{p}(f_{1},\ldots,f_{n},f_{n+1})}\,\Big)$$
et la preuve du th\'eor\`eme est achev\'ee.
\end{proof}

Le th\'eor\`eme \ref{teor:central} est un cas particulier du r\'esultat 
\`a suivre.

%%%%%%%%%%%%%%%%%%%%%%%%%%%%%%%%%
\vskip10pt
\begin{thm}%\label{}
\'Etant donn\'e deux familles \,$\mathcal{F},\mathcal{G}\subset \mathcal{V}$ 
des \difeos telles que le sous-groupe de \,$\mathrm{Diff}^{1}(\RR)$\, engendr\'e par 
\,$\mathcal{F}\cup\mathcal{G}$\, soit ab\'elien et supposons qu'il existe un point 
\,$p\in \mathrm{Fix}(\mathcal{G})$\, dont la 
\,$\mathcal{F}$-orbite soit born\'ee. Alors, il existe un point
\,$q\in\mathrm{Fix}(\mathcal{F},\mathcal{G})$\, dans 
\,$\mathrm{Conv}\big(\,\overline{\mathcal{O}_{p}(\mathcal{F})}\,\big)$\,.

\end{thm}
\vskip5pt
%%%%%%%%%%%%%%%%%%%%%%%%%%%%%%%%%

\begin{proof}[Preuve]
On commence par rappeler l'identit\'e
$$\text{Fix}(\mathcal{G})\cap\text{Fix}(\mathcal{F}) \cap \text{Conv}\Big(\,\overline{\mathcal{O}_{p}(\mathcal{F})}\,\Big) = 
\bigcap_{\substack{f\in\mathcal{F} \\ g\in \mathcal{G}}}\Big\{
\text{Fix}(g) \cap \text{Fix}(f) \cap 
\text{Conv}\Big(\,\overline{\mathcal{O}_{p}(\mathcal{F})}\,\Big) \Big\}
$$ 
o\`u
\,$\text{Fix}(g)\cap\text{Fix}(f) \cap 
\text{Conv}\Big(\,\overline{\mathcal{O}_{p}(\mathcal{F})}\,\Big)$\, est compacte.
Le th\'eor\`eme \ref{teo:cas:fini} nous garantit que 
$$\text{Fix}(\mathcal{G'})\cap \text{Fix}(\mathcal{F'}) \cap \text{Conv}\Big(\,\overline{\mathcal{O}_{p}(\mathcal{F})}\,\Big) \neq\emptyset $$
lorsque les sous-familles \,$\mathcal{F'}\subset\mathcal{F}$\, et 
\,$\mathcal{G'}\subset\mathcal{G}$\, sont finies. Cela d\'emontre que la famille
des compacts 
$$\Big\{\text{Fix}(g)\cap\text{Fix}(f) \cap 
\text{Conv}\Big(\,\overline{\mathcal{O}_{p}(\mathcal{F})}\,
\Big)\Big\}_{f\in\mathcal{F},\;g\in\mathcal{G}}$$
a la propri\'et\'e de l'intersection finie. Par cons\'{e}quent,
$$\text{Fix}(\mathcal{G})  \cap \text{Fix}(\mathcal{F}) \cap \text{Conv}\Big(\,\overline{\mathcal{O}_{p}(\mathcal{F})}\,\Big) \neq\emptyset $$
ce qui donne le r\'esultat recherch\'e.
\end{proof}

%%%%%%%%%%%%%%%%%%%%%
\vskip30pt
\section{Preuve du Lemme Topologique}
\vskip10pt
%%%%%%%%%%%%%%%%%%%%%

La preuve du lemme topologique se fait en deux \'etapes.

\vskip10pt

%:Etapa I
\noindent
{\it \'Etape I\,.}

Supposons tout d'abord que \,$\Gamma=\Gamma(a_{1},\ldots,a_{n})$\, a seulement un nombre fini 
de points d'auto-intersections. Dans ce cas, quitte \`a augmenter de fa\c con convenable le nombre de sommets de la courbe, on peut admettre que 
si deux segments qui composent la courbe s'intersectent, l'intersection se fait sur les extr\'{e}mit\'{e}s. Ainsi,
chaque point de la courbe  \,$\Gamma$\, est recouvert une seule fois, sauf les sommets 
qui peuvent l'\^etre plusieurs fois.

On consid\`{e}re maintenant un point \,$b$\, de  \,$\Gamma$\, o\`u on a 
des intersections. On a vu dans la section \ref{lema:top} comment modifier la courbe de telle fa\c con \`a d\'{e}truire toutes les intersections persistantes par de petites perturbations qui peuvent arriver au point \,$b$. Apr\`es avoir fait une premi\`ere modification,  la courbe initiale peut avoir d'autres points avec des 
auto-intersections persistantes non encore d\'etruites.

Comme on a pr\'{e}serv\'{e} l'orientation initiale des segments,  l'indice de chaque point  \,$q\in\RR-\Gamma$\, par rapport \`a la courbe modifi\'ee co\"{\i}ncide avec celui de la courbe initiale \,$\Gamma$.

R\'ep\'etons le processus de d\'{e}formation de la courbe  initiale, un nombre fini de fois, jusqu'a \'{e}liminer tous les auto-inter\-sec\-tions persistantes par de petites perturbations de \,$\Gamma$. On aura  construit \`a la fin un nombre fini de courbes simples ferm\'ees 
\,$\gamma_{1},\ldots,\gamma_{m}$\, satisfaisant les conditions 
$(i)$\,, $(ii)$\,, $(iii)$ et $(iv)$. Les justifications sont les suivantes.

Pendant le processus de d\'{e}formation, par construction,  les courbes restent toujours des courbes
 ferm\'ees du type \,$\Gamma(b_{1},\ldots,b_{s})$, contenue dans \,$\Gamma$\, et chaque point du plan est recouvert le m\^eme nombre de fois par 
 \,$\Gamma$ et par \,$\Gamma$ modifi\'e; cela signifie que les conditions
 \,$(i)$\, et \,$(ii)$\, sont satisfaites. 
 
 De plus,  l'angle
entre deux segments quelconques de la courbe initiale qui s'intersectent est inf\'{e}rieur \`a
\,$\pi/2$. Donc, si une des courbes \,$\gamma_{1},\ldots,\gamma_{m}$\, poss\`ede au moins une auto-intersection non persistante par de petites perturbations, alors on  d\'{e}duit qu'il existe encore des intersections 
persistantes  car la courbe est ferm\'ee et on est dans le plan. Cela contredit notre hypoth\`ese sur la non existence d'intersections persistantes. On illustre cette situation dans les deux figures  suivantes.

\vskip10pt

\setlength{\unitlength}{1mm}%
\hfil
%\fbox{%
\begin{picture}(48,40)(-1.5,-3)%

\put(0,0){\includegraphics[width=4cm]{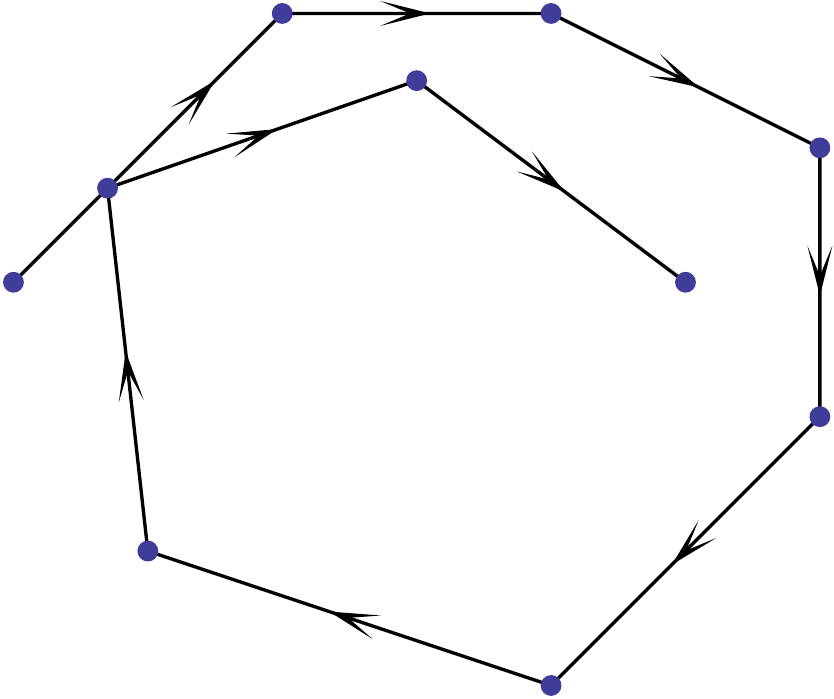}}

\put(-1,17){\scriptsize$b_{\nu}$}
\put(-2,25){\scriptsize$b_{\nu+1}$}
\put(12,34.5){\scriptsize$b_{\nu+2}$}
\put(27,33.5){\scriptsize$b_{\nu+3}$}
\put(40,26){\scriptsize$b_{\nu+4}$}
\put(40,13){\scriptsize$b_{\nu+5}$}
\put(27,-2){\scriptsize$b_{\nu+6}$}
\put(3,5){\scriptsize$b_{\nu+7}$}
\put(6,22){\scriptsize$b_{\nu+8}$}
\put(21,30){\scriptsize$b_{\nu+9}$}
\put(29,18){\scriptsize$b_{\nu+10}$}

\end{picture}%
%}
\hfil
%\fbox{%
\begin{picture}(47,40)(0,-3)%

\put(0,0){\includegraphics[width=4cm]{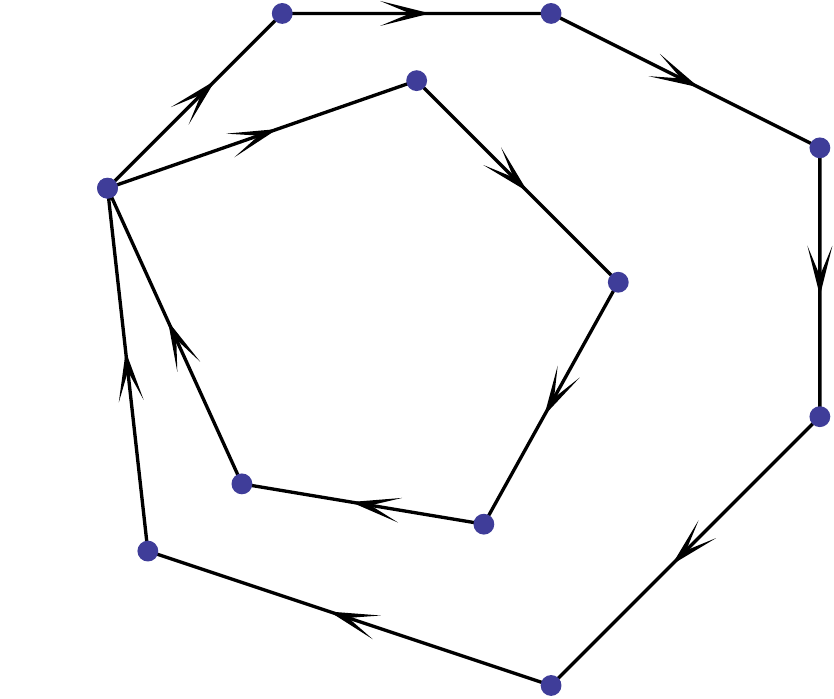}}

\put(12,34.5){\scriptsize$b_{\nu+1}$}
\put(27,33.5){\scriptsize$b_{\nu+2}$}
\put(40,26){\scriptsize$b_{\nu+3}$}
\put(40,13){\scriptsize$b_{\nu+4}$}
\put(27,-2){\scriptsize$b_{\nu+5}$}
\put(3,5){\scriptsize$b_{\nu+6}$}
\put(1,23){\scriptsize$b_{\nu}=b_{\nu+7}=b_{\nu+12}$}
\put(21,30){\scriptsize$b_{\nu+8}$}
\put(30,18){\scriptsize$b_{\nu+9}$}
\put(20,5){\scriptsize$b_{\nu+10}$}
\put(12,11){\scriptsize$b_{\nu+11}$}

\end{picture}%
%}

\vskip10pt

Si \,$\gamma_{i}$\, et \,$\gamma_{j}$\, s'intersectent avec \,$i\neq j$\, alors les int\'{e}rieurs des disques bord\'es par ces courbes ou bien sont disjoints, ou bien l'un d'entre eux est contenu dans l'autre, car il n'existe pas 
d'intersections persistantes. Cela d\'emontre que la famille \,$\gamma_{1},\ldots,\gamma_{m}$\,
a la propri\'et\'e $(iii)$. 

Pour la propri\'et\'e $(iv)$, on a d\'ej\`a remarqu\'e qu'au cours de chaque modi\-fi\-ca\-tion  r\'ealis\'ee sur \,$\Gamma$,
l'indice de chaque point \,$q\in \Gamma-\RR$\, n'a pas chang\'e.

Il reste \`a prouver que la propri\'et\'e $(v)$ est satisfaite.

Pour cela, on consid\`{e}re une \,{\it famille maximale}\, de disques ferm\'es dans 
\,$\RR$
$$D_{n_{1}}\subsetneq D_{n_{2}}\subsetneq \cdots \subsetneq D_{n_{k}}$$
dont les bords sont des courbes de la famille 
\,$\gamma_{1},\ldots,\gamma_{m}$\,. Une telle famille de disques est maximale dans le sens suivant:
\begin{enumerate}
\item[--] 
il n'existe pas de disque ferm\'e ni proprement contenu dans 
\,$D_{n_{1}}$, ni  contenant  proprement le disque \,$D_{n_{k}}$\,, dont le bord soit une courbe de la famille $\gamma_{1},\ldots,\gamma_{m}$\,;

\item[--]
il n'existe pas de disque ferm\'e  proprement contenu dans 
\,$D_{n_{j+1}}$\, et contenant proprement le disque \,$D_{n_{j}}$\,, dont le bord soit une courbe de la famille $\gamma_{1},\ldots,\gamma_{m}$\, pour
chaque \,$j\in\{1,\ldots,k-1\}$\,. 

\end{enumerate}

\'{E}videmment, pour chaque famille de disques 
ferm\'es et embo\^{\i}t\'{e}s
$$D_{s_{1}}\subsetneq D_{s_{2}}\subsetneq \cdots \subsetneq D_{s_{\ell}}$$
dont les bords sont des courbes
de la famille \,$\gamma_{1},\ldots,\gamma_{m}$\,, il existe toujours une famille maximale de disques ferm\'es qui la contient.

S'il existe une famille maximale qui se r\'{e}duit \`a un seul disque, alors la propri\'et\'e $(v)$\, est vraie. Admettons maintenant que toute famille maximale a 
au moins deux disques. Supposons, par contradiction, que 
\begin{align}\label{ind:prov:nnul}
\sum_{j=1}^{k}\text{Ind}_{q}(\partial D_{n_{j}})=0 \quad 
\text{pour tout} \quad q\in \text{Int}(D_{n_{1}})
\end{align}
et pour toute famille maximale de disques ferm\'es
\,$D_{n_{1}}\subsetneq D_{n_{2}}\subsetneq \cdots \subsetneq D_{n_{k}}$\,.

On sait que \,$\text{Ind}_{q}(\partial D_{n_{k}})\neq0$\, pour tout 
\,$q\in \text{Int}(D_{n_{1}})$\,. Admettons, sans perte de g\'{e}n\'{e}ralit\'{e}, que cet
indice est positif et prenons \,$1\leq\mu<k$\, le plus grand entier tel que 
\,$\text{Ind}_{q}(\partial D_{n_{\mu}}) <0$\,. Un tel entier existe  \`a cause
de l'\'{e}galit\'{e} \eqref{ind:prov:nnul} et parce que
\,$\text{Ind}_{q}(\partial D_{n_{k}})>0$\, pour tout 
\,$q\in \text{Int}(D_{n_{1}})$\,.

La courbe \,$\Gamma$\, est connexe, donc il existe une courbe \,$\gamma_{\nu}$\,  de la famille
\,$\gamma_{1},\ldots,\gamma_{m}$\, contenue dans l'anneau
\,$D_{n_{\mu+1}}-\text{Int}(D_{n_{\mu}})$\,,  qui touche 
\,$\partial D_{n_{\mu}}$\, et qui est distincte de 
\,$\partial D_{n_{\mu}}$\, et \,$\partial D_{n_{\mu+1}}$\,.  
Comme l'angle entre deux segments de la courbe initiale qui s'intersectent est inf\'{e}rieur \`a \,$\pi/2$\,,  les courbes 
\,$\partial D_{n_{\mu}}$\, et \,$\partial D_{n_{\mu+1}}$\, ne s'intersectent  pas pour une raison d'orientation. 
D'autre part, la maximalit\'e de la famille de disques nous permet de d\'{e}duire 
que la courbe \,$\gamma_{\nu}$\, borde un disque \,$D_{\nu}$\, contenu dans
\,$D_{n_{\mu+1}}-\text{Int}(D_{n_{\mu}})$\,. 
De plus, on aura aussi que
\,$\partial D_{\nu}$\, a l'orientation anti-horaire, comme \,$\partial D_{n_{\mu+1}}$\,.

Maintenant, on consid\`{e}re une famille maximale de disques ferm\'es qui contient 
\,$D_{\nu}\,,D_{n_{\mu+1}}$\, et \,$D_{n_{k}}$:
\begin{align}\label{fam:disq:02}
\cdots \subsetneq D_{\nu}\subsetneq\cdots\subsetneq D_{n_{\mu+1}}\subset \cdots \subset D_{n_{k}}\,.
\end{align}

Dans cette famille, on peut avoir des disques qui contiennent 
proprement \,$D_{\nu}$\, et qui sont contenus proprement dans 
\,$D_{n_{\mu+1}}$\,, mais les bords de ces disques sont orient\'es dans le sens
anti-horaire vu que tous ces disques doivent toucher le disque 
\,$D_{n_{\mu}}$\,.

\vskip10pt

\setlength{\unitlength}{1mm}%
\hfil
%\fbox{%
\begin{picture}(40,33)(0,0)%

\put(0,0){\includegraphics[width=4cm]{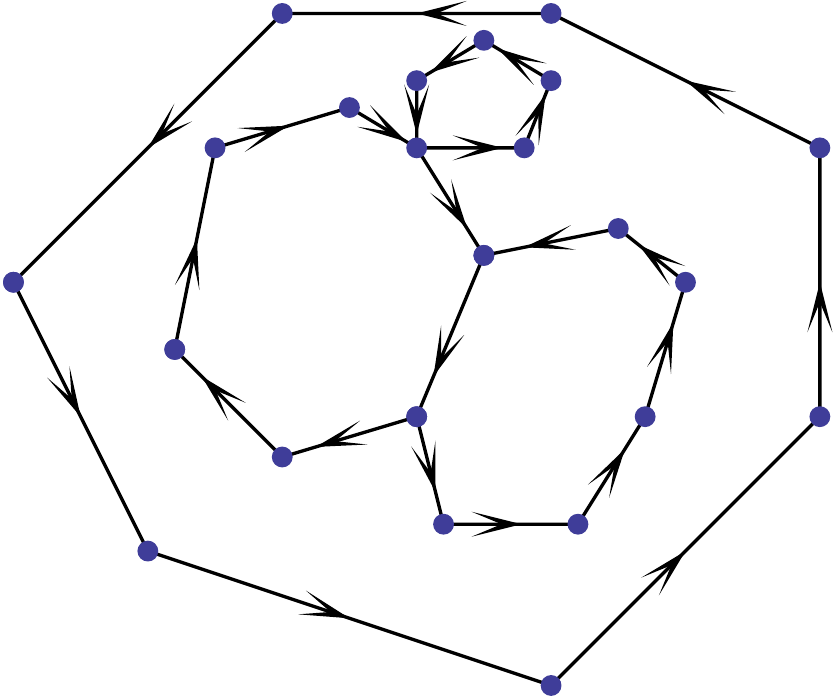}}

\put(30,26){\scriptsize$D_{n_{\mu+1}}$}
\put(12,20){\scriptsize$D_{n_{\mu}}$}
\put(23,14){\scriptsize$D_{\nu}$}
\put(33,17){\scriptsize$\gamma_{\nu}$}
\end{picture}%
%}

\vskip10pt

Maintenant, on r\'{e}p\`{e}te pour la famille \eqref{fam:disq:02} le m\^eme raisonnement utilis\'e
pour la premi\`ere famille maximale de disques. De cette fa\c con, la condition
\eqref{ind:prov:nnul}  pour les familles maximales de disques nous montrera, apr\`es un nombre fini de pas, que la famille \,$\gamma_{1},\ldots,\gamma_{m}$\, poss\`{e}de plus que \,$m$\, courbes, ce qui nous donne la contradiction finale cherch\'ee et d\'emontre que la famille \,$\gamma_{1},\ldots,\gamma_{m}$\, a la propri\'et\'e $(v)$.

\vskip10pt

%:Etapa II
\noindent
{\it \'Etape II\,.}

Pla\c cons-nous sur le cas g\'{e}n\'{e}ral. Quitte \`a augmenter le nombre de sommets de la courbe, on peut admettre que si deux segments \,$[a_{i},a_{i+1}]$\, et \,$[a_{j},a_{j+1}]$\,
qui composent la courbe initiale s'intersectent, alors nous avons l'une des deux possibilit\'es:
\begin{enumerate}
\item[--]
soit l'intersection se r\'{e}duit a un  point unique et dans ce cas le point est l'extr\'{e}mit\'{e} commun aux deux segments;

\item[--] 
soit l'intersection contient plus d'un point et dans ce cas les segments 
co\"{\i}ncident.

\end{enumerate}

Consid\'{e}rons alors un segment \,$[a\,,b]$\, qui compose la courbe 
\,$\Gamma=\Gamma(a_{1},\ldots,a_{n})$\, et supposons qu'il est recouvert
exactement \,$k\geq2$\, fois par la courbe \,$\Gamma$. Avec une petite perturbation de \,$\Gamma$\,, on peut s\'{e}parer le segment 
\,$[a\,,b]$\, dans \,$2k$\, segments
\,$\big\{[a\,,\lambda_{\ell}],[\lambda_{\ell}\,,b]\big\}_{1\leq\ell\leq k}$\,
de telle fa\c con que les \'{e}l\'{e}ments de la famille 
$$\big\{[a\,,\lambda_{\ell}]\cup[\lambda_{\ell}\,,b]\big\}_{1\leq\ell\leq k}$$
s'intersectent deux \`a deux seulement dans les points \,$a\,,b$\, et que
\,$\big([a\,,\lambda_{\ell}]\cup[\lambda_{\ell}\,,b]\big) \cap \Gamma=\{a,b\}$\, pour tout \,$\ell\in\{1,\ldots,k\}$\,.
 
 \vskip10pt

\setlength{\unitlength}{1mm}%
\hfil
%\fbox{%
\begin{picture}(37,47)(-3,-3)%

\put(0,0){\includegraphics[width=3cm]{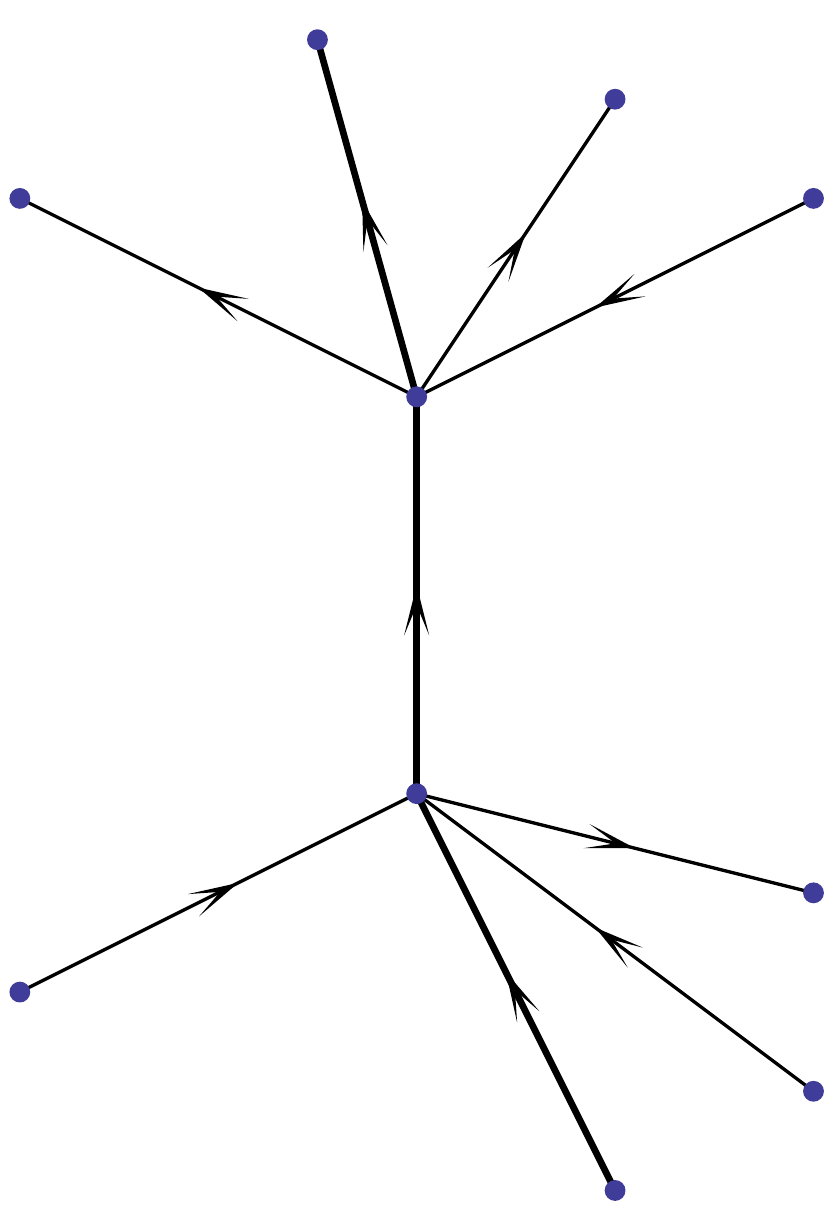}}

\put(3,12){\tiny$_{(1)}$}
\put(16,5){\tiny$_{(2)}$}
\put(24,5){\tiny$_{(1)}$}
\put(24,14){\tiny$_{(1)}$}
\put(16,22){\tiny$_{(3)}$}
\put(25,33){\tiny$_{(1)}$}
\put(20,36){\tiny$_{(1)}$}
\put(13,40){\tiny$_{(2)}$}
\put(3,32){\tiny$_{(1)}$}
\put(12,15){\scriptsize$a$}
\put(12,27){\scriptsize$b$}
%\put(29,18){\scriptsize$b_{\nu+10}$}

\end{picture}%
%}
\hfil
\setlength{\unitlength}{1mm}%
\hfil
%\fbox{%
\begin{picture}(37,47)(-3,-3)%

\put(0,0){\includegraphics[width=3cm]{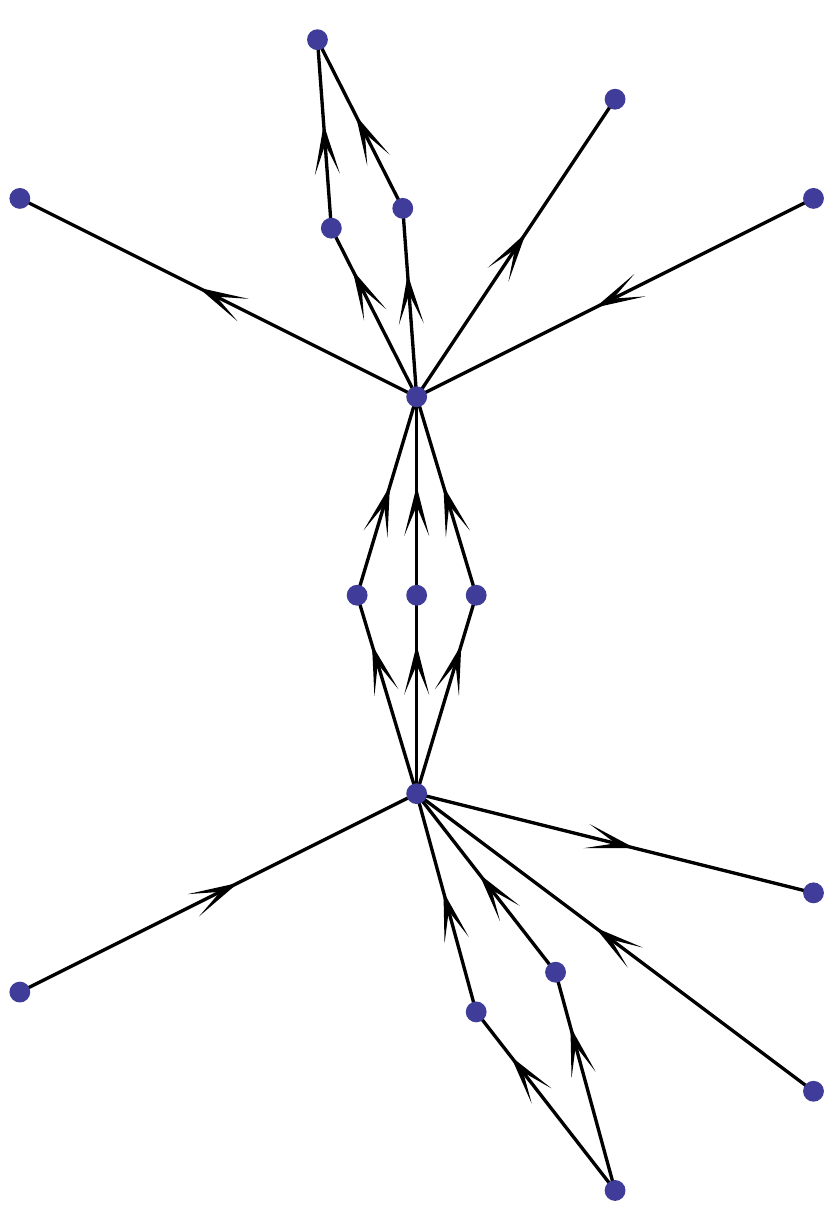}}

\put(12,15){\scriptsize$a$}
\put(12,27){\scriptsize$b$}

\put(9,21.5){\tiny$\lambda_{1}$}
\put(18,21.5){\tiny$\lambda_{3}$}

\end{picture}%
%}

\vskip10pt

Dans cette construction, on suppose aussi que 
\,$[a\,,\lambda_{\ell}]\cup[\lambda_{\ell}\,,b]$\, est toujours \`a gauche de 
\,$[a\,,\lambda_{\ell+1}]\cup[\lambda_{\ell+1}\,,b]$\, o\`u
\,$\ell\in\{1,\ldots,k-1\}$\,. 
De plus, comme on a fait une perturbation petite, on peut supposer que la courbe 
\,$\Gamma$\, ne passe pas par la r\'{e}gion limit\'ee par les segments
\,$[a\,,\lambda_{\ell}],[\lambda_{\ell}\,,b],
 [a\,,\lambda_{\ell+1}],[\lambda_{\ell+1}\,,b] \,.$
 Nous montrons une telle perturbation dans les deux figures ci-dessus. Les nom\-bres \`a c\^{o}t\'{e} de chaque segment indiquent le nombre de fois que le segment est recouvert par la courbe \,$\Gamma$.

On r\'{e}p\`{e}te ce processus avec tous les segments qui sont recouverts plus d'une fois par la courbe initiale. On obtient une courbe 
\,$\Gamma(a_{1}',\ldots,a_{n'}')$\,  dans les m\^emes conditions qu'on avait dans 
l'{\it \'Etape I}\, et alors, on consid\`{e}re la d\'{e}composition de cette courbe 
dans une famille de courbes simples ferm\'ees
\,$\gamma_{1}',\ldots,\gamma_{m'}'$\,
satisfaisant les conditions du lemme.

Finalement, on remarque que chacun des segments 
\,$[a\,,\lambda_{\ell}],[\lambda_{\ell}\,,b]$\, qui appara\^{\i}t dans la perturbation du segment \,$[a\,,b]$\, appartient \`a une unique courbe simple ferm\'ee
\,$\gamma$\, de la d\'{e}composition \,$\gamma_{1}',\ldots,\gamma_{m'}'$\,.
Dans ce cas, on peut faire une isotopie de 
\,$[a\,,\lambda_{\ell}]\cup[\lambda_{\ell}\,,b]$\,
au segment original \,$[a\,,b]$\,, \`a extr\'{e}mit\'{e}s \,$a\,,b$\, fixes
et pr\'{e}servant encore  la propri\'et\'e que la courbe reste toujours simple et ferm\'ee.
De cette fa\c con, on obtient la d\'{e}composition cherch\'ee et la preuve du lemme est termin\'ee.

%\newpage
\vskip30pt
%\bibliography{sapobib} 
%\bibliographystyle{amsplain}

\end{document}